\documentclass{amsart}

\newtheorem{theorem}{Theorem}[section]

\newtheorem{proposition}[theorem]{Proposition}
\newtheorem{nonsense}[theorem]{}
\newtheorem{corollary}[theorem]{Corollary}
\newtheorem{question}[theorem]{Question}

\theoremstyle{definition}

\theoremstyle{remark}

\def \smash {\wedge}
\def \colim {\mathop{\mathrm {colim}}}

\def \L {\mathop{\mathrm {L}}}

\def \C {\mathbb{C}}
\def \F {\mathbb{F}}

\def \rC {\mathrm{C}}
\def \G {\mathrm{G}}
\def \H {\mathrm{H}}
\def \bbH {\mathbb H}

\def \R {\mathbb{R}}
\def \T {\mathbb{T}}
\def \Z {\mathbb{Z}}

\def \CA {{\mathcal A}}

\def \CC {{\mathcal C}}

\def \CE {{\mathcal E}}
\def \CF {{\mathcal F}}
\def \CG {{\mathcal G}}

\def \CP {{\mathcal P}}
\def \CU {{\mathcal U}}
\def \CV {{\mathcal V}}
\def \sc {spin^c}

\def \d {\delta}
\def \f {\phi}

\def \s {\sigma}

\def \Q {\mathbb{Q}}

\def \fs {\mathfrak{s}}

\begin{document}

\title[Refined SW-invariants]{Refined Seiberg-Witten Invariants}
\author[Stefan Bauer]{Stefan Bauer} 
\address{ Fakult\"at f\"ur Mathematik,
Universit\"at Bielefeld, PF 100131, D-33501 Bielefeld}
\curraddr{Institute for Advanced Study, 
Einstein Drive, Princeton, NJ 08540}
\email{bauer@mathematik.uni-bielefeld.de}
\thanks{This material is based upon work supported by the National
Science Foundation under agreement No. DMS-0111298. Any opinions, findings and
conclusions or recommendations expressed in this material are those of the author
and do not necessarily reflect the views of the National Science Foundation.}
\date{\today}
%
\maketitle



In the past two decades, gauge theoretic methods became 
indispensable when considering
manifolds in dimension four. Initially, research centred 
around the moduli spaces of Yang-Mills
instantons. Simon Donaldson had introduced the instanton 
equations into the field. Using
cohomological data of the corresponding moduli spaces, he 
defined invariants
which could effectively distinguish differentiable structures 
on homeomorphic manifolds.
Some years later, Nathan Seiberg and Ed Witten introduced 
the monopole equations.
In a similar spirit as in Donaldson theory, cohomological 
data of the corresponding moduli
spaces went into the definition of Seiberg--Witten invariants. 
These new invariants turned 
out to be far easier to compute, seemingly carrying the same information on 
differentiable structures as Donaldson's. His report \cite{Donaldson96}
gives a glimpse of the wealth of insights in 4-manifold topology 
that could be extracted from these invariants.

However, there is more information in the monopole equations than is seen by the
Seiberg-Witten invariants. The additional information is 
due to an interpretation of the
monopole equations in terms of equivariant stable homotopy. The fact that certain
partial differential equations admit a stable homotopy interpretation is not at all
surprising. Indeed, this has been known for decades \cite{Schwarz}. 
The good news is that in the 
case of the monopole equations it actually is possible 
to make effective use of this fact.
The stable homotopy approach to the monopole equations does not 
only give a different 
view on known results, but also new insights. 

This article is a mixture of a survey and a research article. It serves the
multiple aims of introducing to this area of research, carefully outlining its
foundations, presenting the known results in a unified framework and, last but
not least, proving new results. 

The new results concern various improvements to the definition
of the refined invariants in \cite{BauerFuruta}. 
Theorem \ref{Fredholm}, for example,  
specifies a class of nonlinear Fredholm 
maps between certain infinite dimensional manifolds and shows that the 
path connected components of the space of all such maps are naturally
described by stable cohomotopy groups. This makes it possible to 
define the refined Seiberg-Witten invariant 
as the homotopy class of the monopole map in a 
precise way, clarifying a point left open 
in \cite{BauerFuruta}. The  proof  also indicates how to avoid 
ad hoc arguments used in \cite{Bauer}.
 
Another improvement is on the assumption $b^+>b_1+1$, which had been necessary 
in \cite{BauerFuruta} for a comparison with Seiberg-Witten invariants.
The situation is now summarized in Theorem \ref{Summary}. 
The relation to Seiberg-Witten invariants
is clarified without any restriction on $b^+$ or $b_1$. 
This includes in particular
the wall-crossing phenomenon in the 
$b^+=1$ case, which had been missing in \cite{BauerFuruta}, and the case $b^+=0$.

\bigskip
\section{The monopole map}

The main part in the story to be told is figured by the monopole map
$$
\mu:\CA\to \CC,
$$
which is defined for a closed Riemannian 4-manifold $X$ 
after fixing a $K$-orientation,
or equivalently both an orientation in the usual sense 
and a $\sc$-structure $\fs$. In 
addition, also a background $\sc$-connection has to be fixed. 
The monopole map then is a fiber
preserving map between infinite dimensional vector bundles over the torus
$$
Pic^{\fs}(X)\cong H^1(X;{\R})/H^1(X;\Z).
$$
The refined invariant by definition is the homotopy class of the monopole map in a
sense to be made precise in the next chapter.
This homotopy class does not depend on the chosen 
Riemannian metric or the chosen $\sc$-connection as these choices are parametrized
by connected (indeed contractible) 
spaces and so, indeed, becomes an invariant of the $K$-oriented 
differentiable manifold $X$ (cf. \ref {Universe}). 

Spinors are a main requisite in the definition of the monopole map. Let's start
with the spinor group. 
The group $Spin^c(4)$ consists of those pairs $(u^+,u^-)$ of unitary rank two 
transformations which have the same determinant. If $\Delta^+$ and $\Delta^-$ 
denote the two 
dimensional unitary representations on which the respective
factors act, then the $Spin^c(4)$-representation 
${\mathrm Hom}_{\C}(\Delta^+, \Delta^-)$
admits a real structure. The choice of a basis for the real part $H$ in this 
representation leads to a surjection $Spin^c(4)\to SO(4)$ with
kernel isomorphic to the group $\T$ of complex numbers of unit length.
An element $h$ of $H$ has an adjoint $h^*$ and acts on 
$\Delta=\Delta^+\oplus\Delta^-$ via 
$h(\delta^+, \delta^-)=(-h^*(\delta^-), h(\delta^+))$. This action extends to an
action of the Clifford algebra generated by $H$, resulting in an isomorphism
$Cl(H)\otimes_{\R} {\C}\to {\mathrm End}_{\C}(\Delta)$ of 
$Spin^c(4)$-representations.
Combining with the complexified inverse to the isomorphism 
$Cl(H)\to \Lambda(H)$,
which maps the product $h_1h_2$ to $h_1\smash h_2 - \langle h_1,h_2\rangle$, 
one obtains an isomorphism 
\[
\Lambda_{\C}(H) \to {\mathrm End}_{\C}(\Delta)
\]
of $Spin^c(4)$-representations.  The
decomposition $\Delta=\Delta^+\oplus \Delta^-$ is preserved by  elements 
of   $\Lambda^2_{\C}(H)$.
The kernel of the induced linear map
$$
\rho: \Lambda^2_{\C}(H)\to {\mathrm End}_{\C}(\Delta^+)
$$
consists of the anti-selfdual part $\Lambda^-_{\C}(H)$, its image of the
traceless endomorphisms. The map $\rho$ preserves the real structure, mapping
the real selfdual part $\Lambda^+(H)$ isomorphically to the traceless skew Hermitian
endomorphisms of $\Delta^+$.

We may globalize the above identifications of $Spin^c(4)$-representations to
identifications of bundles by taking fibred products with a principal 
$Spin^c(4)$-bundle. Particularly interesting are such principal bundles
which arise as $Spin^c(4)$-reductions of the othonormal oriented frame bundle 
on an oriented Riemannian four-manifold $X$. These are called 
$\sc$-structures. In fact, the following data do characterize a $\sc$-structure: 
Rank two Hermitian vector bundles $S^+$ and $S^-$, together with
isomorphisms ${\mathrm det}(S^+)\cong {\mathrm det}(S^-)$ and 
$T^*_{\C}X \to {\mathrm Hom}_{\C}(S^+,S^-)$ of Hermitian bundles, the latter 
isomorphism preserving the real structures. In fact, such $\sc$-structures always
exist (compare below). Taking tensor products with 
Hermitian line bundles results in a free and transitive action of $H^2(X;{\Z})$ on
the set of all $\sc$-structures.

There is an interpretation of $\sc$-structures which is special to manifolds
up to dimension 4: Choosing a $\sc$-structure on $X$ is equivalent to choosing
a stably almost complex structure, i.e. an endomorphims $I$ on the Whitney sum
of the tangent bundle with a trivial rank two bundle over $X$ satisfying 
$I^2=-{\mathrm{id}}$. This is because the natural map
$$
BU\to BSpin^c
$$
between the respective classifying spaces induces isomorphisms of homotopy 
groups up to dimension 5. By a theorem of Hirzebruch and Hopf 
\cite{HirzebruchHopf}, there always exists a stably almost complex structure  
on an oriented 4-manifold. If it comes from an
(unstably) almost complex structure, then its second Chern class equals the 
Euler class of $X$. 
Using the equality $c_1^2-2c_2=p_1$ of characteristic classes,
we derive as a necessary condition for a stably almost 
complex structure to be almost complex that its first Chern class satisfies 
$c_1^2= 3\,sign(X) + 2\,e(X)$. If $X$ is connected, this 
condition is also sufficient \cite{HirzebruchHopf}. In this case the integer
\begin{equation}\label{k}
k= \frac{c_1^2 -sign(X)}4 - (b^+ -b_1+1)
\end{equation}
thus measures, how far a stably almost complex structure is
away from being almost complex. Here $b_1$ denotes the first Betti number of $X$ and
$b^+=\frac12 (b_2+sign(X))$ is the dimension of a maximal linear subspace of 
the second de Rham group of $X$ on which the cup product pairing is positive 
definite. 

After fixing a background $\sc$-connection $A$, a  
$\sc$-structure on $X$ allows to define a Dirac operator
$$
D_A: \Gamma(S^+)\to \Gamma(S^-)
$$
mapping positive spinors, i.e. sections of the Hermitian vector bundle $S^+$,
to negative spinors. 
The local model for the symbol of this operator over a point in $X$ is
obtained by identifying the cotangent space with the real part
$H$ of ${\mathrm Hom}_{\C}(\Delta^+, \Delta^-)$. At each point in $X$ this 
symbol is the generator of Bott periodicity, so it provides a $K$-theory
orientation class (compare \cite{AtiyahBottShapiro}) for the manifold $X$.
Indeed, any $K$-theory orientation of $X$ uniquely arises this way.
The Dirac operator is complex elliptic. Its index is given by
\begin{equation}\label{ind}
{\mathrm{ind}}_{\C}(D_A) = \frac{c_1^2 -sign(X)}8.
\end{equation}

Now fix a $\sc$-structure on the $4$-manifold $X$, 
which from now on will be assumed to be connected 
unless explicitely stated differently.
The gauge group ${\mathcal G}=map(X,{\T})$ 
acts on spinors via multiplication with $u:X\to {\T}$, on $\sc$-connections
via addition of $u{\mathrm d}u^{-1}$. The map sending a pair 
$(A,\phi)$ consisting of a $\sc$-connection and a positive spinor to 
$D_A(\phi)$ is equivariant with respect to the gauge group. 
The action of the gauge group on the space of $\sc$-connections is not 
free. However, restriction to the subgroup ${\mathcal G}_0$
consisting of functions which take value $1$ at a chosen point in $X$ results
in a free action. In particular, the based gauge group ${\mathcal G}_0$ acts 
freely on the affine linear
space $A+ i\,{\mathrm {ker}}({\mathrm d})$, 
where $\mathrm d$ denotes the 
de Rham differential on one-forms on $X$, with quotient
\[Pic^{\fs}(X) \cong H^1(X;{\R})/H^1(X;{\Z}).\] 
The gauge group acting trivially on forms, we obtain $\mathcal G$-spaces 
\begin{eqnarray*}
{\widetilde{\CA}}&=&(A+i\,{\mathrm{ker}}({\mathrm d}))\times\left(\Gamma(S^+)\oplus
              H^0(X;{\R})\oplus \Omega^1(X)\right)\\
{\widetilde{\CC}}&=&(A+i\,{\mathrm{ker}}({\mathrm d}))\times
                \left(\Gamma(S^-)\oplus \Omega^0(X)\oplus H^1(X;{\R})\oplus
                \Omega^+(X) \right)
\end{eqnarray*}
consisting of $\sc$-connections, spinors and forms on $X$. 

Consider the map 
$\widetilde\mu:\widetilde\CA\to\widetilde\CC$ defined by
\begin{equation}\label{SWmap}
(A',\,\phi,\,f,\,a)\mapsto 
\left( A',\,D_{A'}\phi+ ia\phi,\,{\mathrm d}^*a+f,\,a_{harm},
\,{\mathrm d}^+a+\s(\phi)\right).
\end{equation}

Here $\s(\phi)$ denotes the trace free endomorphism $i(\phi\otimes\phi^*-
{1\over 2}\|\phi\|^2\,{\mathrm{id}})$ of $S^+$, considered via the
map $\rho$ as a selfdual 2-form on $X$. Restricted to forms, the map is
familiar from Hodge theory: It is linear, injective with cokernel the space
$H^+(X;{\R})$ of harmonic selfdual two-forms on $X$. 
The map $\widetilde\mu$ is equivariant with respect to the action 
of $\CG$. Dividing by the free action of
the pointed gauge group we obtain the monopole map
\[\mu=\widetilde{\mu}/{\mathcal G}_0:{\CA}\to{\CC}\]
as a fiber preserving map between the bundles 
$\CA=\widetilde{\CA}/{\CG}_0$ and $\CC=\widetilde{\CC}/{\CG}_0$
over $Pic^{\fs}(X)$. The preimage of the section
$(A',\,0,\,0,\,0,\,-F^+_{A'})$ of $\CC$, devided by
the residual $\T$-action, is called the moduli space of monopoles.

For a fixed $k > 2$, consider the fiberwise ${\L}^2_k$ Sobolev completion ${\CA}_k$
and the fiberwise ${\L}^2_{k-1}$ Sobolev completion ${\CC}_{k-1}$
of $\CA$ and $\CC$.
The monopole map extends to a continuous map 
${\CA}_k\to {\CC}_{k-1}$ over $Pic^{\fs}(X)$, which will also be denoted by $\mu$.

We will use the following properties of the monopole map:
\begin{nonsense}\label{equivariant}{\rm
It is ${\T}$-equivariant.
}
\end{nonsense}
\begin{nonsense}\label{l+c}{\rm
Fiberwise, it is   
the sum $\mu=l+c$\, of a linear Fredholm 
map $l$ and a nonlinear compact operator $c$.
}
\end{nonsense}
\begin{nonsense}\label{bounded}{\rm
Preimages of bounded sets are bounded.
}
\end{nonsense}

Equivariance is immediate. The action is the residual 
action of the subgroup $\T$ of gauge
transformations which are constant functions on $X$. This group acts by complex
multiplication on the spaces $\Gamma(S^\pm)$ of sections of complex 
vector bundles and trivially on forms. 

Restricted to a fiber, the monopole map is a sum of the 
linear Fredholm operator $l$, consisting of the elliptic operators $D_A$ and
${\mathrm d}^*+{\mathrm d}^+$, complemented by projections to and inclusions of
harmonic forms. The nonlinear part of $\mu$ is built from the bilinear terms 
$a\phi$ and $\s(\f)$. Multiplication  
${\CA}_k\times{\CA}_k\to{\CC}_k$ is continuous for $k>2$. Combined with  
the  compact restriction map ${\CC}_k\to{\CC}_{k-1}$ we gain the claimed compactness
for $c$: Images of bounded sets are contained in compact sets.

Compact perturbations $l+c:\CU'\to \CU$ 
of linear Fredholm maps between Hilbert spaces enjoy a nice topological property:
The restriction to any bounded, closed subset is proper. The argument is 
straightforward: Let $p$ denote a projection to the kernel of $l$. Then 
the restriction of $l+c$ to a closed subset $A\subset \CU'$ factors through
an injective and closed and thus proper map 
$A\to \CU\times \overline{c(A)}\times \overline{p(A)}$,
$a\mapsto (l(a), c(a), p(a))$, a homeomorphism $(u,s,e)\mapsto (u+s,s,e)$ and
the projection to $\CU$, which is proper as the two other factors are compact.

If the bundles $\CA$ and $\CC$ were finite dimensional, then
the boundedness
property would be equivalent to properness. In this infinite dimensional setting,
the argument above can be used the same way as Heine-Borel in the finite
dimensional case to show that the boundedness condition implies properness. 
It turns out that the
ingredients of the compactness proof for the moduli space \cite{Witten}
also prove the stronger boundedness property \cite{BauerFuruta}: The
Weitzenb\"ock formula for the Dirac operator associated to 
the connection $A'=A+ib+ia$ reads
\[
 D^*_{A'} D_{A'}  =  \nabla^*_{A'} \nabla_{A'}  +  {1\over4}s  +  F_{A'}^+.
\]
Applying the Laplacian $\Delta\,|\phi|^2 $ to the spinor part
of an element $(A+ib, \phi ,f , a)$ in the preimage of $\mu$ leads to an estimate    
\[
  \Delta|\phi|^2 \,\,\,\le\,\,\,
  2\,\langle D^*_{A'}D_{A'}\phi-{s\over4}\phi-{1\over2}F_{A'}^+\phi,\,\,\phi\rangle.
\]               
The crucial point is that the term $F_{A'}^+$ can be replaced by an expression
involving $\s(\phi)$ and terms which are straightforward to estimate.
The Laplacian at the maximum is non-negative. Use of this fact
and standard elliptic and Sobolev estimates then lead to an estimate 
$$\|\phi\|_\infty^4 \leq P(\|\phi\|_\infty)$$
with a polynomial $P$ of order $3$. The boundedness property follows easily
from this.

\bigskip
\section{Enter stable homotopy}

In case the first Betti number of $X$ vanishes, the monopole map is a 
map between Hilbert spaces. The boundedness property 
(\ref{bounded}) of $\mu$ is equivalent
to the statement that $\mu$ extends continuously to a map $S^\CA\to S^\CC$ between 
the one-point completions, where the neighbourhoods of the points at infinity are 
the complements of bounded sets. As spaces, these one-point completions are 
infinite dimensional spheres. 
The monopole map thus rightly may be considered as a continuous map between spheres. 

In the general case, we use a trivialisation 
$\CC\cong  Pic^{\fs}(X)\times \CU$ of the bundle $\CC$ to compose the 
monopole map with the projection $p$ to the fiber $\CU$. Now the boundedness
property of $\mu$ translates as follows: The map $p\circ\mu$ extends 
continuously to a map $T\CA \to S^\CU$ from the Thom space of $\CA$ to the
sphere $S^\CU$.

The idea of the refined invariant is to take the homotopy classes of these
one-point completed maps. As it stands, this idea is of course nonsense: All the 
spaces involved are contractible, even equivariantly. So there is no interesting
homotopy theory.

However, not all is lost. Restriction to maps satisfying
not only (\ref{bounded}), but also property (\ref{l+c}) actually does the trick. 
We will consider the situation in a slightly more general setup.

Let $\CE$ and $\CF$ be infinite dimensional 
Hilbert space bundles over a compact base $B$. 
The structure group is the orthogonal group with its norm topology.
Consider the set $\CP (\CE,\CF)$
of fiber-preserving continuous maps $\phi:\CE\to \CF$ satisfying  (\ref{l+c}) and 
(\ref{bounded}).
Let's equip $\CP(\CE,\CF)$ with the topology induced by the metric 
$$
d(\phi,\psi)= \sup_{e\in \CE} \|j\phi(e)-j\psi(e)\|,
$$
where $j:\CF\to \R\times \CF$ denotes the embedding 
$
f\mapsto (1+f^2)^{-1}(1-f^2, 2f)
$
into the unit sphere bundle in $\R \times \CF$ over $B$. (Actually, there are
various topologies on $\CP(\CE,\CF)$ for which the following theorem is true;
the choice made here is just to be definite.)
Choosing a trivialisation $\CF\cong B\times \CU$ of the bundle $\CF$, 
the path components of $\CP(\CE,\CF)$ roughly can be described through a 
bijection
$$
\pi_0(\CP(\CE,\CF))\cong \coprod_{\alpha\in KO(B)}\pi^0_\CU(B;\alpha).
$$
This description uses stable cohomotopy groups of $B$ with
``twisted coefficients''. These groups need some explanation and as it stands,
the statement is rather imprecise. ``For the purposes of planning strategy'' (\cite{Adams84})
it is useful, to think of this decomposition as presented over the group $KO(B)$. For 
the purpose of rigorous definitions and proofs, much  more care has to be taken.

Let's start from the beginning, from pointed spaces.
The prototype of a topological space with a distinguished base point, usually
denoted by $*$, is the one-point compactification $S^U$ of a finite dimensional 
real vector space $U$ with the point at infinity as base point. The smash product 
$A\smash C$ of pointed spaces is the quotient of their product obtained
by identifying $A\times \{*\}\cup \{*\}\times C$ to a point. In this way
$S^U\smash S^V$ is canonically homeomorphic to $S^{U\oplus V}$. The sphere
$S^{{\R}^n}$ is usually denoted by $S^n$. 
The smash product with $S^1$ induces a functor from pointed spaces to pointed 
spaces, called suspension. 

According to Freudenthal's suspension theorem, which holds for finite dimensional
spaces, iterated suspensions eventually induce isomorphisms of sets of pointed
homotopy classes
$$
[S^n\smash A, S^n\smash C]\to [S^{n+1}\smash A, S^{n+1}\smash C].
$$
The notion of a spectrum arose from the desire to define a category in 
which the elements of the resulting abelian group
$$
\mathop{\mathrm{colim}}\limits_{n\to \infty} \,\,[S^n\smash A, S^n\smash C]
$$
appear as homotopy classes of maps between the objects. There are  various 
ways to construct such categories. The situation suggests to use the 
Spanier-Whitehead category indexed by a universe: 
Objects and morphisms in this category are defined through colimit constructions.
The index category consists of 
the finite dimensional linear subspaces of an infinite dimensional
real Hilbert space $\CU$, called universe, with inclusions as morphisms. 
So an object $A$ in the Spanier-Whitehead category 
associates to $U\subset \CU$ a pointed space $A_U$.
To relate these spaces, we use the inclusion $U\subset W$ to identify $W$ with 
$V\oplus U$, where $V$ is the orthogonal complement to $U$ in $W$. The collection
of spaces $A_U$ comes with identifications 
\begin{equation}\label{Spektrum}
\sigma_{U,W}:S^V\smash A_U\to A_W
\end{equation}
satisfying the obvious compatibility condition
\begin{equation}\label{compatibility}
\sigma_{U,W'}=\sigma_{W,W'}\circ(\mathrm{id}_{S^{V'}}\smash \sigma_{U,W})
\end{equation}
for $W'=V'\oplus W\subset \CU$. The morphism set in the Spanier-Whitehead category 
is the colimit
$$
\{A,C\}_\CU = \mathop{\mathrm{colim}}_{U\subset \CU}\, [A_U, C_U].
$$
over the maps
\[
[A_U,C_U]\,\, \stackrel{\mathrm{id}_{S^V}\smash \_}{\longrightarrow}\,\, 
[S^V\smash A_U, S^V\smash C_U]\leftrightarrow 
[A_W,C_W].
\]
The latter identification is induced by the identifications $\sigma^A_{U,W}$ 
and $\sigma^C_{U,W}$.

Every (homotopy) category of spectra
is supposed to contain some variant
of the Spanier-Whitehead category as a full subcategory. 
So it should do no harm to call the 
objects spectra. It should, however, be pointed out that some authors 
reserve the name spectrum to objects in more elaborate categories.

Any pointed space $A$ canonically defines its suspension spectrum, denoted by $A$
as well, by setting $A_U=S^U\smash A$.

To define objects in the Spanier-Whitehead category, it of course suffices
to define them for a cofinal indexing category, as for example
the subcategory of finite dimensional linear subspaces of $\CU$ containing 
a fixed subspace $U$. So associating for a given pointed space $A$ 
to $W=V\oplus U\subset \CU$ the space $S^V\smash A$  
defines a spectrum different from $A$. 
We may safely denote this desuspension by $\Sigma^{-U}A$.
 
Let $p:\CF\cong B\times\CU\to \CU$ be a trivialisation and
suppose $l:\CE\to\CF$ is a continuous, fiberwise linear Fredholm map.
Let $U\subset \CU$ denote a finite dimensional linear subspace such that 
the index of $l$ is represented by the difference $E-\underline{U}$
of finite dimensional vector bundles on $B$. Here $\underline{U}$ denotes
the trivial vector bundle $p^{-1}(U)$
and  $E=l^{-1}(\underline{U})$.  The one-point compactification  $TE$ 
of $E$ is called Thom space of $E$. 
The Thom spectrum is defined as 
$T({\mathrm {ind}}\,l)=\Sigma^{-U}TE$. With this notation,
stable cohomotopy with twisting ${\mathrm {ind}}\,l$ may be defined by
$$
\pi^0_\CU(B;{\mathrm {ind}}\,l):=\{T({\mathrm {ind}}\,l), S^0\}_\CU.
$$

Such twisted cohomotopy groups are a natural habitat for Euler classes of 
vector bundles. To explain this, let 
$F$ be a finite dimensional vector bundle over $B$. Choosing a bundle
isomorphism $E\oplus F\cong \underline{U}$ and a section $\sigma$ of $F$, 
this section and the projection to fibers together define a map 
$\sigma+{\mathrm{id}}_{E}$ extending continuously to one-point compactifications
$TE\to S^U$. This map then represents the
stable cohomotopy Euler class $e(F)\in \pi^0_\CU(B;-F)$. 

The relation to the Euler class of a bundle in a multiplicative cohomology theory
$h$ is as follows \cite{CrabbKnapp}: A Thom class
$u\in h^r(B;-F)=h^r(TF,*)$ corresponds to an $h$-orientation of $F$.
The $h$-theoretic Euler class is defined by $e_h(F)=\sigma^*(u)\in h^r(B)$. A
generator $1\in \tilde h^0(S^0)$ gives rise to the Hurewicz map
$\pi^0(B;-F)\to h^0(B;-F)$, which associates to a stable pointed map
$\phi:T(-F)\to S^0$ the element $h^0(\phi)(1)$. Using the product 
pairing $h^0(B;-F)\times h^r(B;F)\to \tilde h^r(B)$, the $h$-theoretic
Euler class and the stable cohomotopy one are related by
$$e_h(F)=h^0(e(F))(1)\cdot u.$$

To formulate the theorem, let's introduce for a fixed fiberwise linear Fredholm
operator $l:\CE\to\CF$ the subspace $\mathcal{P}_l(\CE,\CF)$
of $\CP(\CE,\CF)$ consisting of elements $\phi$
such that $\phi-l$ is fiberwise compact. 

\begin{theorem}\label{Fredholm}
A projection $p:\CF\cong B\times\CU\to \CU$ induces a natural bijection
$$\pi_0(\CP_l(\CE,\CF))\cong \pi^0_\CU(B;{\mathrm {ind}}\,l).$$
\end{theorem}

The theorem also handles homotopies by applying it to 
the base space $B\times [0,1]$. Note that the restriction maps
$$\pi^0_\CU(B\times [0,1]; {\mathrm{ind}}\,l)\to 
\pi^0_\CU(B\times \{i\}; {\mathrm{ind}}\,l\vert_{B\times\{i\}})$$
are isomorphisms. So if for example 
$\phi=l+c=l'+c'$ are two different presentations
as a sum, then the constant homotopy $\phi=\phi_t=(1-t)(l+c)+t(l'+c')$
can be used to identify $\pi^0_\CU(B; {\mathrm{ind}}\,l)$ with 
$\pi^0_\CU(B; {\mathrm{ind}}\,l')$.
Under this identification, the element associated to the decomposition $l+c$
of $\phi$ is mapped to the element associated to the decomposition $l'+c'$.

{\bf Proof of \ref{Fredholm}.} 
Let's briefly sketch a proof of the theorem: An element 
$\xi\in \pi^0_\CU(B;{\mathrm {ind}}\,l)$ is represented by a virtual bundle $E-\underline{U}$
over $B$, together with a map $TE\to S^U$. It may be necessary to suspend
the given map in order that it can be replaced by a homotopic map 
for which the preimage of the base point consists only of the base point.
In particular, $\xi$ then is represented by a proper  map $E\to \underline{U}$.
The given embedding of $E$ into $\CE$ and an identification of the orthogonal
complements $E^\perp \subset \CE$ and $\underline{U}^\perp \subset \CF$, results
in an element of $\CP_l(\CE,\CF)$.

On the other hand, for a given element $\phi\in\CP_l(\CE,\CF)$, choose a real number 
$R>0$ and an $\epsilon $ with $0<\epsilon<R$. 
The boundedness property (\ref{bounded}) of 
$\phi$ implies that the preimage under $p\phi$ 
of the ball of radius $R$ in $\CU$ is bounded in $\CE$. 
Using compactness of $B$, this bounded preimage is mapped 
by the fiberwise compact operator $p\circ(\phi-l)$ into a compact subset of $\CU$. 
We may cover this image with finitely many $\epsilon$-balls, 
the centers of which generate a finite dimensional vector space $U\subset \CU$.
After possibly enlarging $U$, we can assume that the virtual bundle 
$E - \underline{U}$ with $E=(pl)^{-1}(U)$ represents ${\mathrm{ind}}\,l$. 
The restriction $p\phi|_E$ by construction misses the sphere $S_{R}(U^\perp)$
of radius $R$ in the orthogonal complement of $U\subset \CU$. 
This map 
$p\phi|_E$ extends to the one-point completions to give  a continuous map 
$$
TE\to S^\CU\setminus S_{R}(U^\perp).
$$
Composition with a homotopy inverse to the inclusion 
$S^U\to  S^\CU\setminus S_{R}(U^\perp)$ defines an element of 
$\{T({\mathrm{ind}}\,l), S^0\}$. It remains of course to be checked that the two 
constructions lead to well defined maps between the sets in the theorem which
are inverse to each other. This is straightforward, but a little tedious.
Well-definedness uses the discussion in \cite{BauerFuruta} and in particular 
lemma 2.3 there. The second construction obviously is left inverse to the first.
To show that it is right inverse, one has to construct paths in $\CP_l(\CE,\CF)$
from an arbitrary element to an element, which can be ``projected'' onto
the image of the first construction.
Such  a path is made explicit through the following homotopy
$\phi_t$ which starts from $\phi=\phi_0$ and ends at $\phi_1$.
It is constant on a disk bundle of radius $Q$ in $\CE$, which contains the
preimage of an $R$-disk bundle in $\CF$. Outside it is defined by
\begin{equation}\label{homotopy}
\phi_t(e)=
\left(\frac{|e|}{Q}\right)^t\,
\phi\left(\left(\frac{|e|}{Q}\right)^{-t}e\right).\hfill\qed
\end{equation}
\medskip

The theorem describes the path-connected components of these mapping spaces
in terms of a disjoint union of algebraic objects. So one can hardly expect 
the algebraic structure to be reflected in the world of Fredholm maps 
by natural constructions. In particular,
addition of elements in the respective cohomotopy groups seems to be difficult
to describe in terms of Fredholm maps without making use of the theorem.
However, two aspects are inherent in the Fredholm setup.

{\bf Remarks:}\\
{\bf 2.2.}\label{vanishing} 
A Fredholm map $\phi\in\CP_l(\CE,\CF)$, for which $p\phi$ is not surjective, 
describes the zero element in the stable cohomotopy group associated to its 
linearization $l$. To see this, recall from the previous section 
that $\phi$ is proper. In particular, the image of $p\phi$ is a closed subset in 
$U$. If a point $u\in \CU$ is in the complement of the image, then so is a whole
$\epsilon$-neighbourhood of $u$. Now apply the construction above for some
$R>|u|+ \epsilon$ and replace $\phi$ by the homotopic $\phi_1$.  
We choose $U$ so that it also contains $u$. 
The map $p\phi_1|_E$, followed by the orthogonal projection to $U$ is proper
and by construction misses $u$. So its one-point compactification is null
homotopic.

{\bf 2.3.} The bijection (\ref{Fredholm}) respects products: 
If $\CE_i$ and $\CF_i$ denote Hilbert space
bundles over a base space $B_i$, then taking products of maps results in
a product 
$$ 
\CP_{l_1}(\CE_1, \CF_1)\times \CP_{l_2}(\CE_2, \CF_2)\to 
\CP_l(\CE_1\times \CE_2, \CF_1\times \CF_2)
$$
with $l={\mathrm{pr}}_1^*l_1\times{\mathrm{pr}}_2^*l_2$.
This product structure is reflected in the stable cohomotopy counterpart. 
The natural smash product of objects in the Spanier-Whitehead category is not
defined within one universe, but comes with a change of universes:
Let $A$ and $C$ be objects in the Spanier-Whitehead categories indexed by
universes $\CU$ and $\CV$. The smash product $A\smash C$ then is an object in
the Spanier-Whitehead category indexed by $\CU\oplus \CV$. It is defined by
$$(A\smash C)_{U\oplus V}:= A_U\smash C_{V},$$
whenever both sides make sense, which is at least for 
a cofinal subcategory of the indexing category.

If Fredholm maps $\phi_i\in\CP_{l_i}(\CE_i, \CF_i)$ represent elements
$\xi_i\in \pi^0_{\CU_i}(B_i;{\mathrm {ind}}\,l_i)$, then the product $\phi_1\times \phi_2$
represents the cohomotopy class $\xi_1\smash \xi_2$, which is an element of the
group 
$\pi^0_{\CU_1\oplus \CU_2}(B_1\times B_2;{\mathrm{ind}}\,l)$.


\bigskip
\section{Some equivariant topology}

Using equivariant spaces and maps throughout, the above concepts carry
over to an equivariant setting in a straightforward manner. An appropriate
reference is \cite{Adams84}. 

The action of a compact Lie group $\G$ on a pointed $\G$-space $A$
fixes the distinguished base point. If $A$ and $C$ are
pointed $\G$-spaces, then the smash product $A\smash C$ obtains a $\G$-action 
by restricting the natural $\G\times \G$-action to the diagonal
subgroup. A $\G$-universe $\CU$ is a Hilbert space on which $\G$ acts
via isometries in such a way that an irreducible $\G$-representation, if 
contained in $\CU$, is so with infinite multiplicity. A complete 
$\G$-universe contains all irreducible representations. 

An object of the $\G$-Spanier-Whitehead category indexed by $\CU$ associates
to a finite dimensional representation $U\subset \CU$ a pointed $\G$-space 
$A_U$.
The morphism set is the colimit
$$
\{A,C\}_\CU^\G = \mathop{\mathrm{colim}}_{U\subset \CU}\, [A_U, C_U]^\G
$$
of $\G$-homotopy classes of equivariant maps. This morphism set is
a group if the $\G$-universe $\CU$ contains trivial $\G$-representations.

An equivariant projection
$p:\CF\cong B\times \CU\to \CU$  need not exist.
This may happen for example,
if the $\G$-action on the (unpointed) base space $B$ 
of the $\G$-Hilbert bundle $\CF$ is nontrivial or if 
the fiber over a fixed point in $B$ does not qualify as a universe.
If such a projection $p$ exists, it induces a natural bijection
$$
\pi_0(\CP_l(\CE,\CF)^\G)\cong \pi^0_{\G,\CU}(B;{\mathrm {ind}}\,l)
$$
as before. 

The stable cohomotopy groups in \ref{Fredholm} and in particular 
their equivariant
counterparts are barely known. To get a rough impression, consider
the case of a point $B=\{P\}$. If we choose a universe with a trivial
$\G$-action, then the twist ${\mathrm{ind}}\,l$ 
is characterized by an integer $i$ and the group 
$\pi^0_{\G,\CU}(P;{\mathrm{ind}}\,l)$ 
can be identified with the $i$-th stable stem 
$\pi_i^{st}(S^0)$. On the other extreme, if we choose $\CU$ to  be
a complete universe, then the isomorphism class of ${\mathrm{ind}}\,l$ 
gives an element in
the representation ring $RO(\G)$. In the case where $l$ is an isomorphism,
 $\pi^0_{\G, \CU}(P;0)$ is isomorphic to the
Burnside ring $A(\G)$ (\cite{tD}, II.8.4). If $\G$ is a finite group, $A(G)$
is the Grothendieck ring of finite $\G$-sets with addition given 
by disjoint union and multiplication given by product. A point with
the trivial $\G$-action on it represents  $1$. 

Understanding the group $\pi^0_{\G,\CU}(P;{\mathrm{ind}}\,l)$
for a virtual representation ${\mathrm{ind}}\,l= V-W$ of $G$ 
in any universe $\CU$ whatsoever boils down to understanding
the homotopy classes of $\G$-maps $f:S^V\to S^W$.
Equivariant $K$-theory provides some information 
in case both $V$ and $W$ are
complex representations.  The method is explained in \cite{tD}, II.5:

Let $K_{\G}(B)$ be the Grothendieck group of equivariant complex vector bundles
over the $\G$-space $B$. For a pointed space $B$ as usual $\tilde K_{\G}(B)$
denotes the kernel of the restriction homomorphism $K_{\G}(B)\to K_\G(*)\cong R(\G)$.
If $V$ is a complex $\G$-representation, then $K_\G(V):=\tilde K_\G(S^V)$ is a free
$R(\G)$-module generated by a Bott class $b(V)$ \cite{Atiyah}. 
The image of the Bott class 
$b(U\oplus V)$ under the restriction homomorphism
$\tilde K_\G(S^{U\oplus V})\to K_\G(S^U)$ is  $e_{K_\G}(V)b(U)$.
This defines the Euler class, which was determined by Segal \cite{Segal} to be
the element
\begin{equation}\label{KEulerclass}
e_{K_\G}(V)=\sum_{i=0}^{dim V}(-1)^i \Lambda^i(V)\in R(\G).
\end{equation}
A pointed $\G$-map $f:S^V\to S^W$ induces a $R(\G)$-linear homomorphism in 
$K_\G$-theory. The image of the Bott class of $W$ is a multiple $a_\G(f)b(V)$
of the  Bott class of $V$. The $K_\G$-theory degree $a_\G(f)$ is an element 
of the complex representation ring $R(\G)$. 

To determine the $K_\G$-degree as a character on $\G$, we have to evaluate
it at elements $g\in \G$. Let $\rC$ denote the closure of the subgroup generated
by $g$. Decompose $V=V_\rC\oplus V^\rC$ into the 
$\rC$-fixed point set $V^\rC$ and its
orthogonal complement. The inclusions of the fixed point sets induce a 
commuting diagram
$$
\begin{array}{ccc}
 K_\rC(W)&\stackrel {f^*}\longrightarrow & K_\rC(V)\\
\downarrow&&\downarrow\\
 K_\rC(W^\rC)&\stackrel{f^{\rC*}} \longrightarrow & K_\rC(V^\rC).
\end{array}
$$
The lower map is multiplication by the degree of the map $f^\rC$ as a map
between oriented spheres, if $V^\rC$ and $W^\rC$ both have the same dimension.
Otherwise it is 
zero. This is because $\tilde K(S^{2n})\cong \Z$ 
and  $\pi_{2n-2m}^{st}(S^0)$ is torsion for $n\not= m$. 

Commutativity of the diagram relates Euler classes and degrees by 
\begin{equation}\label{degree}
e_{K_\rC}(W_\rC)d(f^\rC)=a_\rC(f)e_{K_\rC}(V_\rC).
\end{equation}
The representation $V_\rC$ does not contain a trivial summand. So the 
character $e_{K_\rC}(V_\rC)$ does not vanish at $g$. In particular,
the $K_{\G}$-degree can be computed from the 
ordinary degrees of the restrictions of $f$ to fixed points.

Let's apply this concept to a simple example: 
Consider the group $\T$ of complex numbers
of unit length acting on the representation $V$ with character $z\mapsto n+mz$,

\begin{proposition}\label{degree1}
Let $f:S^{2n}\smash S^{\C^m}\to S^{2n}\smash S^{\C^{m+l}}$ be a 
$\T$ equivariant map such that the restricted map on the fixed points has degree
$d\not= 0$. Then  $l\ge 0$ and in case $l=0$ the degree of $f$ nonequivariantly
on the total space is $d$ as well.
\end{proposition}

{\bf Proof:} For $z\not= 1$, the above equation (\ref{degree})
reads as follows:
$$
(1-z)^{m+l}d= a_{\T}(f)(z)(1-z)^m.
$$
The function $d(1-z)^l= a_{\T}(f)(z)$ is a character in $R(\G)$
only if $l\ge 0$. In case $l=0$,  the $K_\T$-degree and hence the ($K$-)degree
equals $d$.\qed

\bigskip
Algebraic topology provides quite heavy machinery for the equivariant 
world. Basic equipment can be found in \cite{Adams84}, \cite{tD}.
Here is a survival kit:

\begin{nonsense}{\rm
An equivariant isometry $\CU\hookrightarrow \CV$ 
induces a change-of-universe morphism 
$\{\,\_\,,\,\_\,\}^G_\CU\to \{\,\_\,,\,\_\,\}^G_\CV$.
It is bijective, if both universes are built from the same irreducible 
representations. 
}
\end{nonsense}

\begin{nonsense} \label{Wirthmueller}{\rm
 A cofiber sequence $A'\to A\to A''$ of pointed $\G$-spaces induces
long exact seqences  
$$
\ldots \leftarrow 
\{S^i\smash A',C\}^\G_\CU\leftarrow\{S^i\smash A,C\}^\G_\CU\leftarrow
\{S^i\smash A'',C\}^\G_\CU\leftarrow \{S^{i+1}\smash A',C\}^\G_\CU\leftarrow \ldots
$$
$$
\ldots \rightarrow 
\{C,S^i\smash A'\}^\G_\CU\rightarrow\{C,S^i\smash A\}^\G_\CU\rightarrow
\{C,S^i\smash A''\}^\G_\CU\rightarrow \{C,S^{i+1}\smash A'\}^\G_\CU\rightarrow 
\ldots
$$
}
\end{nonsense}

\begin{nonsense}{\rm 
Let $\H<\G$ be a subgroup of finite index. Then there are natural bijections
$$\{A, {\mathrm{res}^\G_\H} C\}^\H_{\mathrm{res}^\G_\H \CU}\leftrightarrow
\{(\G\sqcup *)\smash_\H A , C\}^\G_\CU.$$
$$\{ {\mathrm{res}^\G_\H} C, A\}^\H_{\mathrm{res}^\G_\H \CU}\leftrightarrow
\{ C, (\G\sqcup *)\smash_\H A \}^\G_\CU.$$
Here $\G\sqcup *$ denotes the group $\G$ with a disjoint base point and 
 $(\G\sqcup *)\smash_\H A$ is the orbit space of  $(\G\sqcup *)\smash A$
by the action $(h,(g,a))\mapsto (gh^{-1}, ha)$. Here $A$ is
a Spanier-Whitehead spectrum for the group $\H$ and $C$ 
one for the group $\G$. The first adjointness property follows from a corresponding
property on the space level. The second Wirthm\"uller isomorphism can
be found in \cite{tD}, II.6.14. 
}
\end{nonsense}

\begin{nonsense}{\rm 
For spaces with free $\G$-action, the equivariant cohomotopy is naturally
isomorphic to the non-equivariant cohomotopy of the quotient space. 

In more exact diction, this reads:  
Let $A$ be finite dimensional $\G$-space with a free $\G$-action away from 
the base point and let $C$ be a non-equivariant spectrum, indexed by the 
fixed point universe $\CU^\G$ of the universe $\CU$ 
indexing the suspension spectrum of  $A$. 
For such objects, there is a natural bijection
$$
\{A, j^*C\}^\G_\CU \leftrightarrow \{A/\G, C\}_{\CU^\G}.
$$
Here $j^*C$ denotes the spectrum obtained from $C$, 
considered as a spectrum with trivial $\G$-action, by change of universe
$\CU^\G\hookrightarrow \CU$. 
This property is obvious for spaces. The fact that it carries over to
the equivariantly stable world follows from a careful analysis of
the equivariant suspension theorem: In this situation, it suffices
to suspend with trivial representations only to get into the stable range.
}
\end{nonsense}


\bigskip
\section{Topology of the monopole map}

The first statement in the following theorem summarizes the previous discussion.
\begin{theorem}\label{BaFu} {\rm (\cite{BauerFuruta})} 
The monopole map $\mu:\CA\to \CC$ 
defines an element in the equivariant stable cohomotopy group 
         \[\pi^{0}_{{\T},\CU}\left(Pic^{\fs}(X);{\mathrm {ind}}({ D})-
                                       \underline{ H^+(X;\R)}\right).\] 
            For $b^+> {\mathrm {dim}}(Pic^{\fs}(X))+1$, a homology 
           orientation determines a homomorphism 
            of this stable cohomotopy group to $\Z$, which maps
            $[\mu]$ to the integer valued Seiberg-Witten invariant. 
        \end{theorem}

The universe in this statement is explicitly given as the fiber 
$$
\CU=\Gamma(S^-)\oplus \Omega^0(X)\oplus H^1(X;{\R})\oplus
                \Omega^+(X)
$$
of the bundle $\CC$. The index of the linearization of the monopole map consists
of two summands.
The Dirac operator associated to $\fs$ defines a
virtual complex index bundle ${\mathrm {ind}}(D)$ over the Picard torus. 
The second bundle is the trivial bundle with fiber the $b^+$-dimensional
space of self-dual harmonic forms $H^+(X;\R)$. 
An orientation of $Pic^{\fs}(X)\times H^+(X;\R)$ 
is called a homology orientation. 

Let's define the homomorphism of the stable cohomotopy group to $\Z$.
An element of the stable cohomotopy group is represented by an equivariant
map $f:TE\to S^U$ from the Thom space of a bundle $E$ over $Pic^{\fs}(X)$,
where $E-\underline{U}={\mathrm{ind}}\,l$ is the index of the linearization
$l$ of $\mu$. Let $Ci$ denote the mapping cone of the inclusion $i: TE^{{\T}}\to TE$
of the $\T$-fixed point set.
In the long exact sequence associated to the cofiber sequence
$  TE^{{\T}}\to TE\to Ci$, 
\begin{equation}\label{tai chi}
\pi^{-1}_{\T,\CU}(\Sigma^{-U} (TE^{\T})) 
\to 
\pi^0_{\T,\CU}(\Sigma^{-U} Ci) 
\to 
\pi^0_{\T,\CU}(\Sigma^{-U} TE) 
\to 
\pi^0_{\T,\CU}(\Sigma^{-U} (TE^\T)) 
\end{equation}
the first and last term are vanishing because by assumption 
the dimension of the space $S^1\smash TE^\T$ is less than 
the dimension of the $\T$-fixed point sphere in $S^U$, the difference in dimension
being $b^+- b_1 -1$.  So the map $\mu $ can be described by
a cohomotopy element of $\Sigma^{-U}Ci$. 
The Hurewicz image $h(\mu)$ of this element 
in equivariant Borel-cohomology lies in the relative group
$ H^0_{{\T}}(\Sigma^{-U}TE,\Sigma^{-U}TE^{{\T}})$.
The ${\T}$-action on the pair of spaces $(TE,TE^{{\T}})$
is relatively free. So its equivariant cohomology group
identifies with the singular cohomology $H^*(TE/{{\T}},TE^{{\T}})$ 
of the quotient. After replacing
$TE^{{\T}}$ by a tubular neighbourhood, this is the singular cohomology
of a connected manifold relative to its boundary. 
An orientation of $Pic^{\fs}(X)$ together with the standard orientation 
of complex vector bundles defines an orientation class $[TE]_\T$ in the top
cohomology of this manifold.
Similarly, the chosen homology orientation of $X$ and the 
orientation of $Pic^{\fs}(X)$ determine the orientation of $U$ and thus a 
generator $\Sigma^{-U}[TE]_\T$ of the graded cohomology group
$ H^*_{{\T}}(\Sigma^{-U}TE,\Sigma^{-U}TE^{{\T}})$ in its top grading $*=k$.
This cohomology group is a graded module over the polynomial ring  
$H^*_{{\T}}(*)\cong{\Z}[t]$ in one variable $t$ of degree $2$.  
The homomorphism sought for is zero if $k$ is odd or negative. Otherwise
$t^{\frac{k}2}h(\mu)$ is a multiple of the generator  $\Sigma^{-U}[TE]_\T$.
This multiplicity is the Seiberg-Witten invariant. \hfill \qed

To see what happens in the cases $b^+\leq b_1+1$ not covered by this 
theorem we have to take a closer look at the monopole map and distinguish
different cases.

\subsection{The case $b^+=0$.} 
The choice of a point $P\in Pic^\fs(X)$ induces a restriction map 
$$
\pi^0_{\T,\CU}(Pic^\fs(X);\mathrm{ind}\,l)\to \pi^0_{\T,\CU}(P;\mathrm{ind}\,l)\cong
\{S^{\mathrm{ind}(D)}, S^{H^+(X;\R)}\}^\T_\CU
$$
The index of the Dirac operator $\mathrm{ind}(D)$ consists of
$d=\frac18(c^2-\mathrm{sign}(X))$ copies of the tautological complex 
$\T$-representation. 
The restriction to the $\T$-fixed point set of an element in this group is
an element in the stable stem $\pi^{st}_{-b^+}(S^0)$, which is
trivial except in the case $b^+=0$. In this case the restriction of the
monopole map is a linear isomorphism on the fixed point set. Here is an
immediate consequence, well known from Seiberg-Witten theory:

\begin{proposition}
Let $X$ be an oriented $4$-manifold with $b^+=0$. Then the first 
Chern class of any $K$-orientation on $X$ satisfies $c^2\leq {\mathrm{sign}(X)}$.
\end{proposition}

Otherwise the monopole map represented an element in 
$\{S^{\C^d}, S^0\}^\T_\CU$ for some $d>0$ 
which is of degree $1$ on the $\T$-fixed point set. The existence of such an
element contradicts \ref{degree1}.\hfill\qed

Applying Elkies' theorem \cite{Elkies}, we obtain as a corollary Donaldson's 
theorem:

\begin{theorem}\label{diagonal}{\rm (}\cite{Donaldson83}, \cite{Donaldson87}{\rm )}
Let $X$ be a closed oriented four-manifold with negative 
definite intersection form. Then the 
intersection pairing on $H_2(X;{\Z})/Torsion$ is diagonal.\qed
\end{theorem}

\subsection{The case $b^+=1, b_1=0$.}
In the case $b^+=1$, the Seiberg-Witten invariants depend in a well understood
manner (\cite{Witten}, \cite{LiLiu}) on the Riemannian metric and on an additional 
perturbation parameter. 
To understand the phenomenon, let's illustrate it in a characteristic
example. This example describes the situation in the case of an almost complex 
manifold with $b_1(X)=0$, cf. \cite{BauerFuruta}: 

View the spinning globe as a two-sphere with an ${\T}$-action and choose
the north pole as a base point. As a target space, take  a one-sphere with
trivial action and choose two points on this one-sphere as ``poles'', the 
north pole again as base point. Based equivariant maps from the spinning globe
to the one-sphere are determined by their restriction to a latitude, 
which as an arc is a contractible space. So there is only the 
trivial homotopy class of equivariant such maps. 

In contrast, consider equivariant maps, which take north and south pole to
north and south pole, respectively. The monopole maps for 
all choices of metrics and background connections actually are 
of this type.
Such a map basically wraps a 
latitude $n+{1\over 2}$ times around the one-sphere. 
Choosing a generic point in the one-sphere, the oriented count of 
preimages in a fixed latitude defines in a natural way a map of the set of
relative homotopy classes to the integers. This oriented count, 
however, depends
on the choice of the generic point. It changes by $\pm 1$, if 
the generic point is chosen in the ``other half'' of the one-sphere.

There are two ways to deal with the problem. If one prefers to have homotopy
classes, one may consider the monopole map up to equivariant
homotopy relative to the fixed point. The monopole map then is a well defined
element in the set of all such homotopy classes. However, this set will not be 
a group anymore. There is a comparison map to the integers depending
upon the choice of a ``chamber''. 
An alternative is described below.

\subsection{The case $b^+>1, b_1\not= 0$.}
The restriction of
the monopole map to the $\T$-fixed point set
$\CA^\T=Pic^{\fs}(X)\times (H^0(X;\R)\oplus \Omega^1(X))$
is a product map. On the 
factor $Pic^{\fs}(X)$ it is the identity, on the second factor it is
a linear embedding with cokernel $H^+(X;\R)$.
To take this into account, we would like to construct an equivariant spectrum $Q$
encoding this information. Ideally, this spectrum would be obtained by a
push-out of two maps. The one map describes the inclusion $\CA^\T\to\CA$ of
the fixed point set, the other the projection $\CA^\T \to H^0(X;\R)\oplus \Omega^1(X)$.
Such a push-out seems not available in the category we are working in. Let's try
to define a substitute. 
Suppose $E-\underline{U}$ for $U\subset \Gamma(S^-)$
represents the index of the 
Dirac operator as a virtual complex bundle over $Pic^{\fs}(X)$. Then
let $TE/Pic^\fs (X)$ denote the quotient of the Thom space $TE$, where the 
subspace $Pic^{\fs}(X)$, the image of the zero section in $E$, is identified
to a point. Alternatively, $TE/Pic^\fs (X)$ is described as the unreduced 
suspension of the unit sphere bundle in $E$. 
As a $\T$-space it has two fixed points.
The spectrum  $Q(X,\fs,U)$ then is defined by  
$$Q(X,\fs, U)=\Sigma^{-U-H^+(X;\R)}\left(TE/Pic^\fs(X)\right).$$

Using this spectrum, we obtain a straightforward sharpening of \ref{BaFu}:
\begin{proposition}\label{BaFuQ} 
For sufficiently large $U\subset \Gamma(S^-)$, the monopole map $\mu:\CA\to \CC$ 
defines an element in the equivariant stable cohomotopy group 
         $\pi^{0}_{{\T},\CU}(Q(X,\fs,U)).$ 
            For $b^+> 1$, a homology 
           orientation determines a homomorphism 
            of this stable cohomotopy group to $\Z$, which maps
            $[\mu]$ to the integer valued Seiberg-Witten invariant. 
        \end{proposition}
The proof is a slight variation of that of \ref{BaFu}.
To construct the homomorphism to $\Z$,  use the cofiber sequence
$S^0\to TE/Pic^\fs (X)\to Ci$ of spaces. 
The outer terms in the analogue 
of sequence (\ref{tai chi}) now are vanishing for $b^+> 1$ 
for dimension reasons.\qed

One can prove that for $U$ big enough the groups $\pi^{0}_{{\T},\CU}(Q(X,\fs,U))$
become isomorphic. The description, however, still doesn't look satisfactory.

\subsection{The case $b^+=1, b_1\not=0$}
Consider the cofiber sequence 
$S^0\to TE/Pic^\fs (X)\to Ci$ in the proof of \ref{BaFuQ} and set 
$W=U+H^+(X;\R)\subset \CU$. This  
leads to the analogue of the exact sequence (\ref{tai chi})
\begin{equation*}
\pi^{-1}_{\T,\CU}( \Sigma^{-W}S^0)
 \to 
\pi^0_{\T,\CU}( \Sigma^{-W}Ci)
\to 
\pi^0_{\T,\CU}(Q(U)) 
\to 
\pi^0_{\T,\CU}( \Sigma^{-W}S^0)
\end{equation*}
The last term in this sequence vanishes, but the first term is isomorphic
to $\Z$. 
As in the proof of \ref{BaFuQ}, the Seiberg-Witten construction describes
a homomorphism $\pi^0_{\T,\CU}( \Sigma^{-W}Ci) \to \Z$. 
The choice
of a ``chamber'' in computing the Seiberg-Witten invariant amounts to 
the choice of a null homotopy of the restriction
 $S^0\to S^{H^+(X;\R)}$ of the monopole map to the fixed point set.
Such a null homotopy gives rise to 
a lift of the class of the monopole map to an element in 
$\pi^0_{\T,\CU}( \Sigma^{-W}Ci) $.
The wall-crossing formulas mentioned above can be understood as
describing the degree of the composite map
$$\Z\cong
\pi^{-1}_{\T,\CU}( \Sigma^{-W}S^0)
\to 
\pi^0_{\T,\CU}( \Sigma^{-W}Ci)
\to \Z.$$

\subsection{Summary} To summarize the preceeding discussion, let
$\pi^0_{\T,\CU}(Q,Q^\T)$ denote the group 
$\colim_{U\subset \Gamma(S^-)}\pi^0_{\T,\CU}( \Sigma^{-W}Ci)$.
Note that the spectrum $\Sigma^{-W}Ci$ depends on the chosen 
presentation  $E-\underline{U}$ for the virtual index bundle over $Pic^\fs(X)$.
However, the group above by construction is 
independent of the chosen linear subspace $U\subset \CU$.
The groups $\pi^0_{\T,\CU}(Q(U))$ for big enough $U$ become isomorphic, but
not in a natural way. When writing $\pi^0_{\T,\CU}(Q)$, we tacidly
fix some large $U\subset \Gamma(S^-)$.

\begin{theorem}\label{Summary}
The monopole map $\mu:\CA\to \CC$ for an oriented $4$-manifold $X$ with 
$\sc$-structure $\fs$
defines an element in the equivariant stable cohomotopy group 
         $\pi^{0}_{{\T},\CU}(Q(X,\fs)),$ 
which fits into an exact sequence
$$\pi_1^{st}(S^{H^+(X;\R)})
\stackrel{\alpha}{\longrightarrow} 
\pi^0_{\T,\CU}(Q(X,\fs),Q(X,\fs)^\T)
\stackrel{\beta}{\longrightarrow} 
\pi^{0}_{{\T},\CU}(Q(X,\fs))
\stackrel{\gamma}{\longrightarrow} 
\pi_0^{st}(S^{H^+(X;\R)}).$$        
The Seiberg-Witten homomorphism
$h: \pi^0_{\T,\CU}(Q(X,\fs),Q(X,\fs)^\T)\to \Z$ is determined by the 
choice of a homology orientation and relates the monopole class to 
the integer valued Seiberg-Witten invariant  $h\beta^{-1}([\mu])$
in case $b^+>1$.
For $b^+=1$, the choice of a chamber determines a lift
${[\mu]}^{rel}\in \beta^{-1}([\mu])$ and
$h({[\mu]}^{rel})$ is the corresponding Seiberg-Witten invariant. 
The degree of the map $h\alpha$ describes the effect of wall-crossing
on the Seiberg-Witten invariant.
In case $b^+=0$, finally, $\gamma([\mu])=1$.
\end{theorem}


\bigskip
\section{K\"ahler, symplectic and almost complex manifolds}

The current knowledge about differentiable structures on four-dimensional manifolds
builds on the fact that the gauge theoretic invariants are closely related to
the Cauchy-Riemann equations. Witten explained how
in the case of K\"ahler surfaces Seiberg-Witten invariants can be 
determined by complex analytic methods.
Taubes modified the arguments for the case of symplectic manifolds. 
Various mathematicians consequently studied Seiberg-Witten invariants 
for K\"ahler and 
symplectic manifolds. Cutting-and-pasting methods were developped to transfer
these computations to other almost complex manifolds. These efforts resulted in 
a diverse and fascinating picture. 

The refined invariants have little to add
to this direction in four-manifold theory. 
This section intends to explain why.
For the sake of brevity, let's focus on 
central aspects and let's assume $b^+>1$ in this section. 
As noted in the first section, a 
$\sc$-structure is equivalent to a stably almost complex structure on the
tangent bundle of a four-manifold. In particular, an almost complex manifold
comes with a canonical $\sc$-structure ${\mathfrak s}_{can}$. Any other
$\sc$-structure on the underlying oriented $4$-manifold 
is of the form ${\mathfrak s}_{can}\otimes L$ for some
$L\in H^2(X;\Z)$, represented by a line bundle on $X$.
With this convention the first Chern class
of $\mathfrak s_{can}$ is minus the first Chern class $K_X$ of the cotangent bundle.

\begin{theorem}\label{TaubesWitten}
{\rm (}\cite{Witten}, \cite{Taubes94}{\rm )}
Let $X$ be a  symplectic four-manifold with $b^+>1$. The Seiberg-Witten
invariant for the canonical $\sc$-structure ${\mathfrak s}_{can}$ is $\pm 1$. 
Furthermore, Serre-duality holds in the following form:
$$SW({\mathfrak s}_{can}\otimes L)=\pm SW({\mathfrak s}_{can}\otimes (K_X-L)).$$
\end{theorem}

\begin{theorem}\label{TaubesWitten2}
{\rm (}\cite{Witten}, \cite{Taubes96}{\rm )}
Let $X$ be a  symplectic four-manifold with $b^+>1$.
If for some $L\in H^2(X;\Z)$ the Seiberg-Witten invariant of 
${\mathfrak s}_{can}\otimes L$ is nonvanishing, then
this $\sc$-structure corresponds to an almost complex structure.
\end{theorem}

Witten and Taubes actually prove more than is stated in these theorems:
The monopole map is not surjective, unless there is a pseudo-holomorphic
curve in $X$ which is Poincar\'e dual to the class $L$. The result 
follows by the application of adjunction inequalities \cite{MorganSzaboTaubes}.
By remark (2.2), we get as an immediate consequence:

\begin{corollary}\label{symplectic}
Let $X$ be an oriented four-manifold with $b^+>1$, 
which admits a symplectic structure. 
If the stable cohomotopy
invariant $[\mu]\in \pi^{0}_{{\T},\CU}(Q)$  in {\rm \ref{BaFuQ}} is 
non-vanishing for some $\sc$-structure $\mathfrak s$ on $X$, then $\mathfrak s$ 
describes an almost complex structure on $X$.
\end{corollary}

The refined invariants, when applied to symplectic manifolds,
carry exactly the same information as the Seiberg-Witten invariants. 
This is a consequence of \ref{symplectic} and the following statement.

\begin{proposition}\label{almost complex}
Let $X$ be an almost complex four-manifold with $b^+>1$. Then the homomorphism
$\pi^{0}_{{\T},\CU}(Q)\to \Z$ in {\rm \ref{BaFuQ}} comparing the Seiberg-Witten 
invariant with its refinement is an isomorphism.
\end{proposition}

{\bf Proof.} 
For an almost complex 4-manifold, the ``virtual dimension of the moduli space''
$k$ is zero (\ref {k}). The construction of the comparison homomorphism
in \ref{BaFu} and  \ref{BaFuQ} considers a map from a pair 
$(TE/\T,TE^\T)$ of spaces to a sphere.  The integer $k$ is exactly the difference
of the dimensions of $TE/\T$ and the sphere. The dimensions being equal and
$(TE/\T,TE^\T)$ being a connected and oriented manifold relative to its boundary,
one can apply a classical theorem of Hopf. It states that the homotopy classes
of such maps are classified by their degree.\qed

So in order to test, whether the refined invariants are of any use, we have to leave
the by now familiar world of symplectic or at least almost complex 4-manifolds
and enter the jungle.


\bigskip
\section{Some stable cohomotopy groups}

The groups $\pi^{0}_{{\T},\CU}(Pic^\fs(X);{\mathrm {ind}\,l})$ 
seem to be at least as hard to compute
as the stable homotopy groups of spheres. Let's restrict to the simplest cases.
In particular, let's only consider
 $4$-manifolds $X$ with vanishing first Betti number and $b^+>1$. 
The groups then are then determined by the index of the linearization $l$
of the monopole map. We will write $\pi^0_{\T,\CU}({\mathrm {ind}\,l})$
for short.
The index of the Dirac operator is denoted by  
$d={\mathrm {ind}}_\C(D)= \frac{c^2-\mathrm {sign}(X)}8.$ The virtual
dimension (\ref{k}) of the moduli space  is $k=2d-b^+-1$.
\begin{proposition}\label{proj} {\rm (}\cite{BauerFuruta}{\rm )}
              Let $X$ be a $K$-oriented, closed
             $4$-manifold with vanishing first Betti number and $b^+>1$. 
            The stable equivariant cohomotopy group $\pi_{{\T},\CU}^0({\mathrm {ind}\,l})$
             is isomorphic to the nonequivariant stable cohomotopy group
             $\pi^{b^+-1}( P(\C^{d}))$ of the complex $(d-1)$-dimensional
             projective space.
             This group vanishes for $k<0$.
             It is isomorphic to $\Z\oplus A(k,d)$, if $k\ge0$ is even,
             and to $A(k,d)$ otherwise. Here $A(k,d)$ denotes a finite abelian
             group. For any prime $p$, the $p$-primary part of $A(k,d)$ vanishes
             for $k< 2p-3$.
For $k\leq 4$, the groups $A(k,d)$ can be described 
as follows: \begin{itemize}
\item $A(0,d)\cong A(4,d)=0$. 
\item $A(1,d)\cong A(2,d)$. For even $d$ these groups are isomorphic to $\Z/2$, 
otherwise they vanish. 
\item The $2$-primary part of $A(3,d)$ is a cyclic group, depending on the congruence
class of $d$ modulo $8$. The order of the group is $8,0,2,4,4,0,2,2$ for 
the congruence classes $0,1,2,\ldots$.
\item The $3$-primary part of $A(3,d)$ is of order $3$ if $d$ is divisible by $3$ and
else vanishes.
\end{itemize}
\end{proposition}     
The proof of the first statement uses the sequence (\ref{tai chi}), which in this
situation by excision is a part of the long exact cohomotopy sequence
for the pair $(D({\C}^d)\sqcup *, S({\C}^d)\sqcup *)$
consisting of the unit ball and and its bounding 
sphere in the complex vector space ${\C}^d$ with an extra base point added.
The $\T$-action on the sphere is free, so we may apply (3.5) to get the result.
The Atiyah-Hirzebruch spectral sequence accounts for the rest of the statement.

Instead of chasing through technicalities, let's try to understand in an informal
way, how to represent elements in these groups for small $k$. 
Recall the structure of the stable homotopy groups
of spheres in low dimensions. The group $\pi^{st}_n(S^0)$ is cyclic for $n\le 5$.
It is infinite for $n=0$, of order $2$ for $n=1$ or $2$, 
of order $24$ for $n=3$ and zero else in this range. For $n=1$ and $3$, these groups
are generated by Hopf maps $S(\F^2)\to P(\F^2) $ for $\F =\C$ and $\F=\bbH$,
denoted by $\eta$ and $\nu$. These generators satisfy the relation $\eta^3=12\nu$.

First consider the map obtained by forgetting the
$\T$-action. This homomorphism 
$$f:\pi_{\T,\CU}^0({\mathrm {ind}\,l})\to\pi_{k+1}(S^0).$$
associates to a $\T$-equivariant map between $\T$-representation spheres its
underlying nonequivariant map. In the case $k=0, d=2$, the group 
$$\pi_{\T,\CU}^0({\mathrm {ind}\,l})\cong \{S^{\C^2},S^3\}^\T_\CU\cong \pi^2( P(\C^2))\cong \Z$$ 
is generated by the unreduced suspension of the Hopf map $\eta$.
 For $k=0$ and general $d$, the
generator of $\pi_{\T,\CU}(Q)\cong \Z$ is mapped to $(d-1)\eta$.

The collapsing map $P(\C^{d})\to P(\C^{d})/P(\C^{d-1})
\cong S^{2d-2}$ induces a homomorphism
$$c:\pi_k(S^0)= \pi^{b^+-1}(S^{2d-2}) \to \pi^{b^+-1}(P({\C}^{d}))\cong
\pi_{\T,\CU}({\mathrm {ind}\,l}).$$
This map turns out to be an isomorphism for $k=0$ and surjective onto the
torsion subgroup $A(k,d)$ for $0<k\leq 4$. 
The composite map $f\circ c:\pi_k(S^0)\to \pi_{k+1}(S^0)$ 
is multiplication by $(d-1)\eta$. 

Finally consider the Hurewicz map 
$$\pi^{b^+-1}(P(\C^{d}))\to H^{b^+-1}(P(\C^{d})).$$
If an element in $\pi^{b^+-1}(P(\C^{d}))$ represents a monopole 
map, then the image of this element under the Hurewicz map is a multiple of
the generator in the cohomology group in the respective dimension. 
This multiplicity is the Seiberg-Witten invariant. 

The Hurewicz map is neither surjective nor injective, the kernel being torsion.
The non-injectivity issue makes the stable cohomotopy invariant a true refinement
of Seiberg-Witten invariants. This will be addressed in the next sections.
Non-surjectivity implies that, depending on $k$ and $d$, the Seiberg-Witten
invariants automatically satisfy certain divisibility conditions. 
The index of the image of the Hurewicz map
$\pi^{2m}(P(\C^{m+n}))\to H^{2m}(P(\C^{m+n}))$
for $m,n\ge 0$ is known to be the stable James number $U(-m,n)$ 
(cf. \cite{CrabbKnapp88}, Remark 2.7). These James numbers can be 
defined in a more general setup and appear in various geometric situations. 
$K$-theory methods provide an estimate for them,
which conjecturally is sharp:
\begin{theorem} {\rm (}\cite{CrabbKnapp88}{\rm )}
The power series in $z$ with rational coefficients
$$\left(\frac{z}{\log(1+z)}\right)^m,$$
when multiplied with $U(m,n)$, becomes integral modulo $z^n$. 
\end{theorem}


\bigskip
\section{Intermezzo}

This chapter aims at sensitizing for some snags one should be 
aware of when working in this field. 
One concerns a misinterpretation of the Pontrijagin-Thom construction, 
another the proper use of homotopy categories. 

The main difference between the
familiar approach to gauge theory and the homotopy approach is the replacement
of spaces by maps. The Pontrijagin-Thom construction provides a perfect and well-known 
duality between the concepts ``stable homotopy classes of maps 
between spheres'' and  ``bordism classes of framed manifolds''.
At first glance, this duality suggests stable maps to contain equivalent information
as localized data in the form of  moduli spaces 
together with suitably specified normal bundle data. 
This idea is particularly appealing to anybody
working in gauge theory, since the use of localized data 
--often in form of characteristic classes--
is a main trick of the trade.
I propose to dispose of this idea as quickly as possible, since it
is prone to deception and self-deception. 
Here is a much too long discussion; for related topics compare \cite{Adams84}, ch. 6.

\subsection{Equivariant transversality.}
One minor reason is due to the fact that the Pontrijagin-Thom correspondence fails
in general in an equivariant setting due to the fact that transversality
arguments don't work in sufficient generality. 

\subsection{A variant on the Eilenberg-swindle.}\label{Eilenberg}
More seriously, the information cannot possibly localize as suggested above. 
The reason is as puzzling as it is simple: Any two framings on a bundle
by their very definition are isomorphic.
Framings can be distinguished only embedded in a surrounding space. 

But how to keep control over framings when changing the surrounding space?
The default surrounding space we are dealing with is a 
Hilbert space. In order to get into business, we have to reduce to finite dimensions.
And, to get this straight, the only natural way is by linear projection. 
Indeed, such projections are used in the proof of \ref{Fredholm}. Now comes the point: 
Embedded framed manifolds are extremely ill behaved under projections. 

Let's look at this in more detail. Embed $S^1$ as the 
$\T$-orbit of a nonzero element in $\C^d$ and fix a framing of an affine normal disk
to a given point in $S^1$. Using the $\T$-action on $\C^d$, this framing extends
to a framing of the normal bundle of $S^1$. 
By equivariant (here it is okay) Pontrijagin-Thom, this framing corresponds 
to a generator in the corresponding 
equivariant stable homotopy group, which happens to be isomorphic to $\Z$, as we have seen
in the preceeding chapter.
Now consider a generic projection $\C^d\to \C^{d-1}$. This is a $\T$-equivariant 
map and the $\T$-equivariant normal framing in $d$ complex dimensions is
equivariantly projected to one constructed the same way in 
complex $d-1$ dimensions, which also represents a generator in the corresponding 
group. The disastrous effect on the framing becomes apparent only after
forgetting the $\T$-action.  As explained in the preceeding chapter,  
nonequivariantly the constructed framing of the embedding in $\C^d$ 
is $(d-1)\eta\in \pi_1^{st}(S^0)\cong (\Z/2)\eta$. So it is
trivial for odd $d$ and nontrivial for even. 
In particular, when projecting along an infinite dimensional Hilbert space
in an uncontrolled manner, we systematically do Eilenberg-swindles.
There are several ways to deal with this. I'll explain some commonly used ones.

\subsection{Equivariance to the rescue.}\label{Equivariance}
The way to gain control is by the use of the stable map representing the
framing. Let's do that.
This is an equivariant map $S^{\C^d}\to S^{2d-1}$. From the equivariant picture 
it is clear that this map has nothing to do with equivariant maps
$S^{\C^{d-1}}\to S^{2d-3}$:
Projection should correspond to desuspension. 
But considering source and target, we immediately realize that if our example were
a desuspension, then it were along different $\T$-representations on either
side of the map. The lesson should be that only
by holding to the map as a double-entry book-keeping device, we can tell
legal and harmless projections (desuspensions) from the illegal and harmful.
But actually, in our case this is not enough.

\subsection{Universes to the rescue.} 
Let's take a closer look at the example just discussed and let's forget that there 
was a $\T$-action. As pointed out, linear projections should 
correspond to desuspensions. But if we forget the $\T$-action, 
the linear projections in the example on both sides are real linear along an $\R^2$.
As we have seen, they cannot correspond to desuspensions. 
Intuitively, the problem is easy to understand:
In the source, we are trying to desuspend a ``moving frame'', whereas in the
target, we want to desuspend a ``fixed frame''. 
Now that we have excluded representation theory
to act as a savior, we need a replacement to convey that idea. 
The notion of a universe, which seems to
go back to Peter May, is such a replacement. The point here is that 
the projection above along $\R^2$ does not factor through a projection along $\R^1$ as it 
should. If one uses universes, this feature is built in.

\subsection{On the usage of spectra I}\label{spectraI}
I want to present a way how not to define the refined invariants: 
This uses the spectrum of a self-adjoint elliptic operator, acting on a Hilbert space $\CU$.
After choosing an oriented basis for eigenspaces,
we get a canonical embedding $\R^\infty\to\CU$, which we may use to make suspensions
ordered by the integers instead of finite dimensional linear subspaces of $\CU$. This
is okay if one does not change the operator. 

The snag appears if one wants to change the operator.
Let's do that, say by changing a metric used to define it. 
At first sight this looks controllable: A small 
change of the operator will result in a small change of the eigenvalues, so
locally, up to ``canonical'' homotopy, this should define a ``canonical''
homotopy equivalence between the sphere spectra indexed by the integers. 

Will this stand up to scrutiny?
Assume we have a closed path
of operators such that the eigenspaces for eigenvalues in a fixed interval
constitute a bundle over $B=S^1$. It may happen that the bundle for a chosen set
of eigenvalues is not orientable. 
Following the ``canonical'' homotopy equivalences of the sphere 
along the circle, we obtain that the
identity map over the base point is ``canonically'' homotopic to minus the identity
map. This is not what we want.

But, the space of metrics is contractible. So
we may always extend the operator to an operator
parametrized by a disk. In the critical cases this will 
involve other eigenspaces than the ones
we started with. So only very special arrangements of eigenspaces will be ``admissible''
for the argument. And which arrangements are ``admissible'' may depend heavily on the
chosen extension of the operator on the disk. There may exist no ``admissible'' arrangement
that works in all situations.

Orientation is governed by a 
determinant line bundle, which exists in the Fredholm setting.
So, indeed, there  may be a way to coherently enforce all such bundles over $B$
to be orientable. I don't know any, but let's suppose we found one. 
Then, as bundles over $B$, they are trivial. However, there
are two trivializations up to homotopy to choose from. 
If we pick the wrong one, we will have the following phenomenon:  
Using the trivialization, we may parallel transport an embedded 
$S^1$ with framed normal bundle in the fiber over a point in $B$ once around 
the loop $B$. This parallel transport changes the framing.

But, the space of metrics is contractible.
Indeed, if the operator is such that the bundle over $B$ extends to a bundle of eigenspaces
over the disc, then this would pick a trivialization.  
However, there may be a different extension of the operator to the disk such that we get a 
trivialization only if we add a $2$-dimensional eigenspace. The two trivializations
obtained that way need not be the same, as the example (\ref{Eilenberg}) shows.
Which to choose?

Let's stop here. Who ever desires to use eigenspaces of self-adjoint 
operators to define homotopy objects is kindly asked not to rely on 
bluff and belief, but on
reason. I cannot see, how to create well-defined mathematical objects this way 
and I doubt it is possible.  I recommend universes instead. 

\subsection{Why do universes work?}\label{Universe} 
The discussion above shows well-definedness of the refined invariants 
to be a non-trivial issue. 
The following argument is not based on the 
contractibility of some parameter space, but
on the contractibility of the orthogonal group of Hilbert space.
If we take a path in our parameter space (metrics, $\sc$-connections),
then we will get a bundle of universes over that path. The theorem
of Kuiper \cite{Kuiper} shows that this is a trivial bundle and has a unique 
trivialization up to homotopy. 
A trivialization identifies the universes defined for different parameters.
Such an identification of universes provides for a change-of-universe 
isomorphism of the stable cohomotopy groups defined with respect to the respective
universes. A trivialization homotopic modulo end points to the one chosen will induce the same
isomorphism of stable cohomotopy groups. This uses \ref{Fredholm} for the parameter space
$B\times [0,1]\times [0,1]$.
In particular, by taking closed paths, we get
that the invariants are well defined.

\subsection{On the usage of spectra II}\label{htop}
Finally, I want to point out a reliable avenue to create nonsense.
Is it possible to construct (homotopy types of) spectra 
out of spaces, which themselves are only defined up to homotopy?
That means, all spaces are defined up to homotopy, the suspensions are defined up 
to homotopy, the compatibility condition (\ref{compatibility})
holds only up to homotopy.
The answer in general is: No. This
would amount to a lift from the homotopy 
category of topological spaces to the category of topological spaces. 
This problem has been addressed
in work of Dwyer and Kan, compare e.g. \cite{DwyerKan1}, \cite{DwyerKan2}.
To see the problems, just assume 
for each $U\subset \CU$, the space $A_U$ to be a sphere homotopy equivalent to $S^U$.
When trying to prove well-definedness of the identity map, not only similar
problems as above turn up, but also higher dimensional phenomena. 
There is no magic to cure this problem. 

Since I am using \cite{tD} as a reference, I should point out that his definition
of spectra looks similar to the one I am criticizing. It actually is different: 
The author 
wisely only uses complex representations as suspension coordinates.
Because of the implicit $\T$-equivariance (\ref{Equivariance}) this gets
rid of all the complications I lamented about. Moreover, the author is only interested
in spectra as realizing equivariant homology and cohomology functors on spaces. 
He does not define
a category of spectra and in particular he does not define maps of spectra, thus 
avoiding any discussion about the indicated higher dimensional phenomena.

Not all authors have taken
this problem in homotopy theory serious. 
Sadly enough, it 
renders a considerable part of the 
literature in this subject useless.


\bigskip
\section{Gluing results}

\subsection{Gluing along positive curvature.}
Connected sums of oriented $4$-manifolds have vanishing 
Seiberg-Witten invariants, unless one of the summands has negative definite
intersection form. The same statement holds for Donaldson invariants.
This fits very well with known stability results on simply connected
$4$-manifolds: A theorem of Wall \cite{Wall} states that if any two differentiable
$4$-manifolds are homotopy equivalent, then after taking connected sum with sufficiently
many copies of $S^2\times S^2$, the resulting manifolds will be diffeomorphic. 
In many cases it is known that already one copy is ``sufficiently many''. For 
example, complete intersections or elliptic surfaces  are almost completely 
decomposable \cite{MandelbaumMoishezon}. 
That means, the result of taking connected sum with a single complex
projective plane is diffeomorphic to a connected sum of projective planes, taken with 
both standard and reversed orientation.
Simon Donaldson defined in \cite{Donaldson89} $mod\,2$-polynomial
invariants, which potentially could distinguish different structures on 
connected sums. However, no examples were found.

The stable cohomotopy invariants don't vanish in general for connected sums. This
shows that they are true refinements of Seiberg-Witten invariants. 
The  connected sum theorem \cite{Bauer}
states that for a connected sum $X_0\#X_1$ of 4-manifolds, the stable
equivariant cohomotopy invariant is the smash product of the invariants of 
its summands. It is straightforward to compute explicite examples.

A precise statement of the theorem constitutes already
a major part of its proof. We will discuss a slightly more general setup.
Let $X$ be the disjoint union of a finite number, say $n$, of closed connected 
Riemannian $4$-manifolds $X_i$, each equipped with a $K$-theory orientation. 
Suppose each component
contains a separating neck $N_i\cong Y\times [-L,L] $. So it is
a union  \[X_i=X_i^-\cup X_i^+\] of closed submanifolds with
common boundary $\partial X_i^\pm=Y\times \{0\}$. Here $Y$ denotes a 
$3$-manifold with a fixed Riemannian structure.
The length $2L>2$ of the neck is considered a variable.
For an even permutation $\tau$ of the indices, let $X^\tau$ be the manifold obtained
from $X$ by interchanging the positive parts of its components, that is
\[{X}_i^\tau=X_{i}^-\cup X_{\tau(i)}^+.\]
Next comes the question of whether and how $K$-orientations glue. In order to be
able to glue,  we of course need the following

{\bf Assumption:} {\it
The $K$-orientations on all components $X_i$, when pulled back along the  
inclusion $ Y\times [-L,L]\stackrel \sim\to N_i\hookrightarrow X_i$,  
lead to the same $K$-orientation.} 

This assumption is automatically satisfied in case $Y$ is an integral homology sphere.
In general, 
in order to get a well-defined $K$-orientation on the manifold $X^\tau$, it does not
suffice to fix an isomorphism class, but we also have to
fix identifications. Note that the gauge group $map\,(Y\times [-L,L], \T)$
acts freely on the set of
all such identifications. If the gauge group is connected, any such identification
will give the same $K$-orientation on $X^\tau$. We can enforce connectedness by the

{\bf Assumption:} {\it Let $Y$ have vanishing first Betti number.}

It turns out that we will have to put much stronger assumptions on the geometry of 
$Y$ in order to prove the gluing theorem. So we need not discuss this tricky issue at
this point.
Under these assumptions, a $K$-theory orientation of $X$ uniquely 
induces by gluing one on $X^\tau$.
A main ingredient for the gluing setup is a change of universe isomorphism
$V_Y:\CU\to \CU^\tau$. Its explicit construction uses a smooth path 
\[\psi:[-1,1] \to SO(n)\] 
starting from the unit, i. e. $\psi(-1)= \mathrm{id}$, and ending 
at $\tau$, considered as the permutation matrix 
$(\d_{i, \tau(j)})_{i,j}\in SO(n)$. 
Suppose we are given a bundle over $X$ such that the restrictions over the necks are
identified with a bundle $F$ over $Y\times [-L,L]$. Using these identifications,
the restrictions of the bundle to $X_i^\pm$ glue together to a bundle over
$X^\tau$.  Sections of the given bundle, when restricted over the neck, can 
be viewed as a section of the  bundle $\oplus_ {i=1}^n F$ over $Y\times [-L,L]$.
Consider the path  $\psi$ as rotation of the components of this bundle. Rotating
 via $\psi$ a given section of a bundle over $X$ 
results in a section of the glued bundle over $X^\tau$. 
This gluing construction, applied to forms and spinors on $X$,
defines fiberwise linear bundle isomorphisms $V_Y:\CA\to{\CA}^\tau$ and 
$V_Y:\CC\to\CC^\tau$ of the Hilbert space bundles over a suitably defined
identification $Pic^{\fs}(X)\stackrel{\cong}{\rightarrow} Pic^{\fs}( {X}^\tau)$.
The following theorem is formulated in \cite{Bauer} only for the cohomotopy groups
in \ref{BaFu} and the case $Y=S^3$. 
The proof extends without further changes to the version in \ref{BaFuQ} and to positively
curved manifolds $Y$, i.e. quotients of the sphere. 
\begin{theorem}\label{Gluing} 
Let $Y$ be a manifold with positive Ricci and in particular scalar 
curvature. Then the change of universe isomorphism  
$$V_Y: \pi^{0}_{{\T},\CU}(Pic^\fs(X);{\mathrm {ind}\,l})\to 
\pi^{0}_{{\T},{\CU}^\tau}(Pic^\fs(X^\tau);{\mathrm {ind}\,l^\tau})$$
identifies the monopole classes of $X$ and ${X}^\tau$
for corresponding $K$-theory orientations.
\end{theorem}

The theorem claims the diagram
$$
\begin{array}{ccc}    \CA     &   \stackrel{\mu}\longrightarrow  &    {\CC}\\
       V\downarrow\phantom{V}&                 &   \phantom{V}\downarrow V\\
   {\CA}^\tau& \stackrel{{\mu}^\tau}\longrightarrow& {\CC}^\tau
\end{array}
$$
to commute up to homotopy, i.e. there is a path in $\CP_l(\CA,\CC)$ connecting
the maps $\mu$ and $V^{-1}\mu^\tau V$.  The difference between the two
maps is a compact operator. So the homotopy need only change the 
compact summand
in the monopole map such that 
at any time during the homotopy the boundedness condition (\ref{bounded})
remains satisfied. 
Control is achieved by the use of Weitzenb\"ock formulas for
both the Dirac operator and the covariant derivative. Positivity of scalar and
Ricci curvature, respectively, along the neck 
provide the necessary estimates on the
spinor and form components during the homotopy.  The estimates 
on spinor and forms finally are tuned by neckstretching. So the theorem 
holds for sufficiently large $L$ and hence for any $L>1$.
The proof in \cite{Bauer} actually constructs a path in a
slightly bigger space than $\CP_l(\CA,\CC)$. This can
be avoided by the use of the homotopy (\ref{homotopy}).

To apply this theorem, let's spell out the following elementary observation.

\begin{proposition}
 {\rm (}\cite{Bauer}{\rm )} Let $X$ be the disjoint union of a finite number of 
$K$-oriented $4$-manifolds $X_i$.
Then the Thom spectrum $T({\mathrm {ind}\,l})$ of the index bundle over $Pic^\fs(X)$ 
is the smash product of the corresponding spectra $T({\mathrm {ind}\,l_i})$
associated to the components and the 
stable cohomotopy class of the monopole map of $X$ is
the smash product 
\[[\mu(X,\fs)]=\wedge_{i=1}^n[\mu(X_i,\fs_i]\in
\pi^{0}_{{\T}^n,\, \oplus \CU_i}(Pic^\fs(X);{\mathrm {ind}\,l})\]
of the stable cohomotopy classes associated to the respective components. 
The action of the torus ${\T}^n$ on the sum $\oplus_{i=1}^n \CU_i$ is factorwise.
\end{proposition}

Note that the $\T$-action in \ref{Gluing} on these spectra is the diagonal one.

The proof of the connected sum formula follows from applying this theorem to
the case $Y=S^3$ when $X$ is the disjoint union of a connected sum $X_0\# X_1$
and two copies of the $4$-sphere. The manifold $X^\tau$ then will be the disjoint
union of $X_0$, $X_1$ and one further copy of the $4$-sphere. 
Using (\ref{Summary}) and (\ref{degree1}) it is immediate to recognize
$\mu(S^4)$ as homotopic to the identity map on the sphere spectrum.
In particular, we may identify the monopole class $[\mu(X_0\#X_1,\fs_0\#\fs_1 )]$,
with the monopole class
$$[\mu(X_0\#X_1,\fs_0\#\fs_1)\smash \mu(S^4)\smash \mu(S^4)]=
[\mu(X_0\#X_1,\fs_0\#\fs_1)\smash {\mathrm{id}}_{S^0}\smash {\mathrm{id}}_{S^0}],
$$
via some obvious change-of-universe identifications.
So as a corollary to \ref{Gluing} we obtain the connected sum theorem.

\begin{theorem}\label{connected sum} {\rm (}\cite{Bauer}{\rm )}
 The gluing map $V_{S^3}$ identifies the class 
$
[\mu(X_0\#X_1,\fs_0\#\fs_1 )]
$ of the monopole map
of the connected sum of two $K$-oriented $4$-manifolds with the smash product
$[\mu(X_0,\fs_0)]\smash[\mu(X_1,\fs_1)]$ of the monopole classes of the summands.
\end{theorem}

The gluing theorem applies to a further construction, 
which is discussed in \cite{GompfStipsicz}, p. 411:
Suppose, the $K$-oriented $4$-manifolds $X_0$ and $X_1$ both contain a 
$-2$-curve, i.e. a smoothly embedded $2$-sphere with self-intersection number $-2$.
Cutting out  tubular neighbourhoods of these $-2$-curves, we obtain manifolds with
real projective $3$-space as boundary. Using an orientation reversing diffeomorphism 
of the boundaries, we may glue the manifolds along their boundaries. Let $X_0\#_2 X_1$
denote the resulting manifold. 

The orientation reversing diffeomorphism permutes the two spin-structures on
the real projective $3$-space $P(\R^4)$. One of these two spin-structures extends
as a spin-structure to the tubular neighbourhood of a $-2$-curve. This property
distinguishes the two spin-structures. The other spin structure extends to a
$\sc$-structure on the tubular neighbourhood, the determinant line bundle of which
has degree congruent $2$ mod $4$, when restricted to the $-2$-curve.

Gluing two copies of tubular neighbourhoods of $-2$-curves by the use
of an orientation reversing diffeomorphism of the boundaries, results in a 
manifold $N$. This manifold can also be recognized as the manifold
 $N=\overline{P(\C^3)\# P(\C^3)}$ obtained by reversing the
orientation on the connected sum of two copies of the complex projective plane. 
There are four $\sc$-structures on $N$ for which the monopole map is homotopic 
to the identity map on the sphere spectrum. This again is immediate from
(\ref{Summary}) and (\ref{degree1}). Exactly the same argument as in \ref{connected sum},
with $S^4$ replaced by $N$, thus proves:

\begin{theorem}\label{connected 2-sum}
Let $X_0\#_2X_1$ be the sum of two $4$-manifolds along $-2$-curves with $\sc$-structure
$\fs$. Then there are $\sc$-structures $\fs_0$ and $\fs_1$ on $X_0$ and $X_1$, respectively,
such that one of the associated first Chern classes evaluates at the corresponding
$-2$-curve with $2$, the other with $0$ and $\fs=\fs_0\#_2\fs_1$. 
The gluing map $V_{P(\R^4)}$ identifies the class of the monopole map
$
[\mu(X_0\#_2X_1,\fs )]
$
with 
$[\mu(X_0,\fs_0)]\smash[\mu(X_1,\fs_1)].$
\end{theorem}

Obviously, the range of applications of \ref{Gluing} is rather limited. 
It would be desirable to extend the stable cohomotopy approach in a well-defined manner 
(cf. \ref{htop}) to manifolds with boundary explaining the behaviour
under cutting and pasting.

\subsection{Applications to $4$-manifolds.}

The computations of the stable cohomotopy groups in (\ref{proj}) can now be combined with
known results on Seiberg-Witten invariants (\ref{almost complex}). Most of the 
following statements are immediate.

\begin{theorem} {\rm (}Vanishing results for connected sums, \cite{Bauer}.{\rm )}
Let $X$ be a connected sum of  oriented $4$-manifolds.  Then the refined 
invariants vanish for any $\sc$-strucure on $X$ in 
the following cases:
\begin{enumerate}
\item The refined invariants vanish for any $\sc$-structure on one of the summands.
\item There are two or more summands which are symplectic and have
vanishing first Betti numbers. Furthermore, one symplectic summand $X_0$ satisfies 
$b^+(X_0) \equiv 1\, mod\, 4$.
\item The manifold $X$ has vanishing first Betti number, $b^+(X)\not\equiv
\,4\, mod\,8$ and is a connected sum of $4$ symplectic manifolds.
\item The manifold $X$ has vanishing first Betti number
 and is a connected sum of $5$ symplectic manifolds.
\end{enumerate}
\end{theorem}

The theorem remains true, if one replaces ``symplectic'' by the weaker assumption
``all $\sc$-structures with non-trivial refined invariants are almost complex''.

\begin{theorem} {\rm(}Vanishing results for sums along $-2$-spheres.{\rm )}
Let $X_0$ and $X_1$ be oriented $4$-manifolds containing $-2$-spheres $C_0$ and
$C_1$, respectively.  Then the refined 
invariants vanish for any $\sc$-strucure on the sum $X_0\#_2X_1$ along these spheres
in the following cases:
\begin{enumerate}
\item The refined invariants vanish for any $\sc$-structure on $X_0$.
\item The first Chern class of any $\sc$-structure on $X_i$, for which the refined
invariant is nonvanishing, gives the same number 
modulo $4$ when evaluated on $C_i$ (for both $i=0,1$).
\item\label{semi-definite} 
The first Chern classes of those $\sc$-structures on $X_i$, which have 
nonvanishing refined invariants, span a linear subspace of $H^2(X_i;\Q)$ on 
which the cup-product is positive semi-definite. 
\item Both $X_0$ and $X_1$ can be equipped with the structure of a minimal K\"ahler
surface with $b^+(X_i)>1$.
\end{enumerate}
\end{theorem}

{\bf Proof.} The first two statements are immediate from \ref{connected 2-sum}:
The assumptions imply one of the two factors in the smash product to vanish.
The fourth statement is a special case of the third. 
The proof of \ref{semi-definite} uses the following, well-known fact: Complex
conjugation in a small tubular neighbourhood of a $-2$-sphere extends
to an automorphism of the $4$-manifold $X_i$ which is constant outside a larger tubular
neighbourhood. The effect in second cohomology
is a reflection on the hyperplane perpendicular to the Poincar\'e dual $PD(C_i)$
of the $-2$-curve. If there was a $\sc$-structure with non-vanishing refined invariant, whose
first Chern class is not perpendicular to $PD(C_i)$, then $PD(C_i)$ were a linear
combination of the first Chern classes of this $\sc$-structure and its reflected 
$\sc$-structure. This would contradict the assumption. As a consequence, the second
condition is satisfied, the number modulo $4$ being $0$.\qed

\begin{question}
Is there a minimal symplectic $4$-manifold with $b^+>1$
for which the first Chern classes 
of $\sc$-structures with non-vanishing Seiberg-Witten invariants span
a linear subspace of $H^2(X;\Q)$ which is not positive semidefinite?
\end{question}

Here are some general non-vanishing results. Of course, no manifold can be on both 
a vanishing list as above and a non-vanishing list. This has some non-trivial implications.
Note that the assumptions are met by symplectic manifolds.

\begin{theorem}\label{nonvanishing}{\rm (}Non-vanishing 
results for connected sums, \cite{Bauer}.{\rm)}
There is 
a $\sc$-structure on the oriented $4$-manifold $X$
for which the associated refined invariant is non-vanishing, if 
one of the following holds:
\begin{enumerate}
\item The manifold $X$ is a connected sum $X=X_0\#X_1$ of
a manifold $X_0$, which admits a $\sc$-structure with non-vanishing refined invariant, and
a manifold $X_1$ with  $b^+(X_1)=0$. 
\item The manifold $X$ has vanishing first Betti number and is a connected sum with two
or three summands. For every summand $X_i$ there is an almost complex structure for 
which the integer Seiberg-Witten invariant is odd and $b^+(X_i)\equiv\,3\,mod\,4$.
\item The manifold $X$ is a connected sum with four summands, 
has vanishing first Betti number and $b^+(X)\equiv\,4\,mod\,8$. 
 For every summand $X_i$ there is an almost complex structure for 
which the integer Seiberg-Witten invariant is odd and $b^+(X_i)\equiv\,3\,mod\,4$.
\end{enumerate}
\end{theorem}

{\bf Proof.}
Only the first statement is not discussed in \cite{Bauer}. Using Donaldson's
theorem \ref{diagonal} we can find a $\sc$-structure on $X_1$ such that the 
virtual index bundle of the Dirac operator over $Pic^\fs(X_1)$ has rank $0$. The inclusion
of a point in $Pic^\fs(X_1)$ induces a restriction map 
$\pi^0_{\T, \CU}(Pic^\fs(X_1);{\mathrm {ind}\,l})\to \pi^0_{\T, \CU}(S^0)$. 
The image of the monopole class
is the identity map.\qed

The information retained in the refined invariants of 
connected sums is much more detailed than
these sweeping vanishing and non-vanishing theorems might suggest. To get an impression,
let's consider connected sums of certain elliptic surfaces which  
had been classified \cite{Bauer92} up to diffeomorphism with methods from Donaldson theory. 
Note that in each of the two homeomorphism classes of such elliptic surfaces 
there are infinitely many diffeomorphism classes. 

\begin{corollary} \label{elliptic}{\rm (}\cite{Bauer}{\rm )}
            Suppose the connected sum $\#_{i=1}^4E_i$ of simply connected
            minimal elliptic surfaces of 
            geometric genus one is diffeomorphic to a 
            connected sum $\#_{j=1}^nF_j$ of elliptic surfaces. 
            Then $n=4$ and the $F_j$ and the $E_i$ are 
            diffeomorphic up to permutation.
\end{corollary}

Ishida and LeBrun \cite{IshidaLeBrun1}, \cite{IshidaLeBrun2} pointed out some 
differential geometric applications of the connected sum theorem. In particular
they proved non-existence statements for Einstein metrics on connected sums
of algebraic surfaces.


\section{Additional symmetries}

\subsection{Spin structures}
The case of spin structures was pioneered by Furuta \cite{Furuta}.
The key observation is that for a spin $4$-manifold $X$ the monopole map
is actually $Pin(2)$-equivariant, where $Pin(2)\subset Sp(1)\subset \bbH$ is
the normalizer of the maximal torus $\T\subset \C\subset\bbH$ in $Sp(1)$. This
subgroup is generated by $\T$ and an additional element $j\in \bbH$ satisfying $j^2=-1$ and
$ij+ji=0$.

The group $Spin(4)$ is isomorphic to the product of two copies of $Sp(1)\cong SU(2)$
and embeds as a subgroup in $Spin^c(4)$. So the $Spin^c(4)$-representations  used
in the definition of the monopole map naturally restrict to $Spin(4)$-representations.
Considered this way as $Spin(4)$-representations, $\Delta^+$ and $\Delta^-$ admit
quaternionic structures. The Dirac operator, therefore, is $\bbH$-linear. 

This additional structure is not preserved by the monopole map:
Consider the induced action of $Sp(1)$ on the space of all $Spin(4)$-equivariant
quadratic maps $\Delta^+\to \Lambda^+$. The isotropy group of the term $\sigma$ in
the definition of the monopole map is $\T$. The normalizer of the torus interchanges $\sigma$
and $-\sigma$. This indicates, for which action and which group the monopole map can be made
equivariant.

Taking the $spin$-connection $A$ as the background $\sc$-connection, we can define a 
$Pin(2)$-action on the spaces $\CA$ and $\CC$ used in the definition of the monopole map: 
The group acts via the quaternionic
structure on the sections of the quaternionic bundles $S^+$ and $S^-$. The element
$j$ acts via multiplication by $-1$ on both forms and $\sc$-connections (after identifying
the space of connections with $A+i\Omega^1(X)$. The monopole map
(\ref{SWmap}) indeed is equivariant with respect to this $Pin(2)$-action.

In this setup, our standard universe 
$\CU=\Gamma(S^-)\oplus \Omega^0(X)\oplus H^1(X;\R)\oplus \Omega^+(X)$
will not contain trivial $Pin(2)$-representations. As a consequence, equivariant cohomotopy
$\pi^0_{Pin(2), \CU}(Pic^\fs(X);{\mathrm {ind}\,l})$ 
in general does not carry a group structure; it is just a set which
even may be empty. Indeed, the main result in \cite{Furuta} in effect proves emptiness of this
set in certain cases.

In order to get groups, we may simply enlarge the universe $\CV=\CU \oplus H$ 
by adding an infinite dimensional Hilbert space $H$ with trivial $Pin(2)$-action.
The change-of-universe map
$\pi^0_{Pin(2), \CU}(Pin^\fs(X);{\mathrm {ind}\,l})\to \pi^0_{Pin(2), \CV}(Pin^\fs(X);{\mathrm {ind}\,l})$ 
can be viewed as induced by
smash product with the identity element in $\pi^0_{Pin(2), H}(S^0)$ and 
turns out \cite{DissBirgit} to be injective in case $b_1=0$;
the image is characterized in algebraic terms. In particular, there is no loss of
information by this change-of-universe, but a gain of convenient algebraic structure.

\begin{theorem}{\rm (}\cite{Furuta}{\rm)} Let $X$ be a $spin$ $4$-manifold with
$sign(X)<0$. Then the second Betti number of $X$ satisfies the inequality
$$b_2(X)\ge 2-\frac{10}8 \,sign(X).$$
\end{theorem}

{\bf Proof.} The inclusion $A\hookrightarrow Pic^0(X)$ of the $spin$-connection
induces a restriction map $\pi^0_{Pin(2), \CV}(Pic^\fs(X);{\mathrm {ind}\,l})\to \pi^0_{Pin(2), \CV}({\mathrm{ind}\,l})$.
The index of $l$ is
$\frac{-sign(X)}{16}\bbH - H^+(X;\R)$. In order to apply the $K$-theory degree formula,
we need to complexify these $Pin(2)$-representations. So, consider the square of the monopole map
$$\nu=[\mu(X)]\smash[i\mu(X)]\in \pi^0_{Pin(2), \CV+i\CV}\left(\Sigma^{-H^+(X;\R)\otimes\C}(S^{\bbH^d})\right).$$
The element $j\in Pin(2)$ acts by multiplication with $-1$ on the $Pin(2)$-representation 
$H^+(X;\R)\otimes\C$. We would like  to compute
the $K_{Pin(2)}$-mapping degree $a_{Pin(2)}(\nu)$ of $\nu$ via the formula (\ref{degree})
$$e_{K_{Pin(2)}}(H^+(X;\R)\otimes\C)\cdot d_{Pin(2)}=a_{Pin(2)}(\nu)\cdot e_{K_{Pin(2)}}(\bbH^d),$$
 which takes 
place in the representation ring $R(Pin(2))\cong \Z[\lambda,h]/(\lambda^2-1,\lambda h-h)$. 
Here $\lambda$ stands for the one-dimensional representation on which $j$ acts by multiplication
with $-1$ and $h$ stands for the quaternions, viewed as a $Pin(2)$-representation.
The singular cohomology mapping degree $d_{Pin(2)}$ can be computed by considering for each element
in $Pin(2)$ the cohomology degree on the fixed point spheres of that element. By dimension
reasons, this vanishes except for the conjugates of $j$. It is $1$ for $j$ itself, as by construction
$\nu$ is the identity on the fixed point set. So we get $d_{Pin(2)}=\frac12(1-\lambda)$. The 
 $K_{Pin(2)}$-Euler classes are computed via (\ref{KEulerclass}) to be
$e_{K_{Pin(2)}}(H^+(X;\R)\otimes\C)=(1-\lambda)^{b^+}$ and $e_{K_{Pin(2)}}(\bbH^d)=(2-h)^d$.
The mapping degree formula thus reads
$$\frac12 (1-\lambda)^{b^++1}=a_{Pin(2)}(2-h)^d.$$
In the representation ring, this equality can only be satisfied, if $a_{Pin(2)}$ is of
the form  $a(1-\lambda)$ for some integer $a$ (the character on the left hand side is zero on $\T$!). 
So we are left with the equation
$$2^{b^+-1}(1-\lambda)=a2^d(1-\lambda),$$
which can be satisfied only for $d\leq b^+-1$, or equivalently, $b^+\geq 1-\frac{sign(X)}{8}.$\qed

This theorem can be sharpened a little:

\begin{theorem} {\rm (}\cite{DiplomBirgit},\cite{DissBirgit}{\rm)} 
 Let $X$ be a $spin$ $4$-manifold with
$sign(X)<0$. Then the second Betti number of $X$ satisfies the inequality
$$b_2(X)\ge 2a-\frac{10}8 \,sign(X),$$
with $a=2$, if $sign(X)\equiv 32 \, mod\,64$ and $a=3$, if $|sign(X)|\equiv 48 \, mod\,64$
and $a=1$ else. Moreover, in the case $sign(X)=-64$, one has $b_2(X)\ge 88$. 
\end{theorem}

The two references correspond to two different proofs. The first one relies on results 
of  S. Stolz and M. Crabb
in $\Z/4$-equivariant stable homotopy.
The second one imitates in principle the proof above, using more refined 
$KO_{Pin(2)}$-mapping degrees instead. Furuta in \cite{Furuta2} announced that by the same methods
one can prove $a=3$ if  $sign(X)\equiv 0 \, mod\,64$. 

The following example shows that
these methods cannot be carried through to prove the so-called $\frac{11}8$-conjecture stating
$b_2(X)\ge -\frac{11}8 sign(X)$.

\begin{theorem}\label{minimalspin} {\rm (}\cite{DissBirgit}{\rm)} 
There is an element in $\{S^{\bbH^5},S^{V^{12}}\}^{Pin(2)}_\CU$, 
where $V$ is the real $1$-dimensional
non-trivial $Pin(2)$-representation. 
\end{theorem}

So the lowest rank of a potential counterexample to the $\frac{11}8$-conjecture, at least
according to current knowledge, is $b_2= 104$. 

The connected sum theorem also works in the $Pin(2)$-equivariant setting for $spin$-manifolds.
Taking connected sum of a $spin$-manifold $X$ with $S^2\times S^2$ amounts for the 
$Pin(2)$-equivariant monopole classes to multiplication with the stable cohomotopy 
Euler class $e(V):S^0\hookrightarrow S^V$ of the $Pin(2)$-representation $V$. 

The long exact sequence for the pair of spaces $(D(V)\sqcup *,S(V)\sqcup *)$, together with the
adjunction \ref{Wirthmueller} leads for a $Pin(2)$-spectrum $A$ 
to a long exact Gysin sequence 
$$\ldots \to \pi^{-1}_{\T, \CV}( S^V\smash A)
\to \pi^0_{Pin(2), \CV}(S^V \smash A)
\to \pi^0_{Pin(2), \CV+V}(A)
\to \pi^{0}_{\T,\CV+V}( A)
\to \ldots .
$$
The map in the middle is multiplication with the Euler class $e(V)$. The next map 
restricts the group action. Application to the Thom spectrum $T({\mathrm {ind}\,l})$ of 
the index bundle over $Pic^\fs(X)$ gives as an 
immediate consequence:

\begin{theorem} {\rm (}\cite{DissBirgit}{\rm)}\label{Birgit} 
Suppose the $Pin(2)$-equivariant monopole class of a $spin$ $4$-manifold $X$ with $sign(X)<0$ 
is not divisible by the Euler class $e(V)$. Then the refined Seiberg-Witten invariant
of $X$ is nonzero. In particular, if $X$ has second Betti number $b_2(X)=-\frac{11}8 sign(X)$
for $|sign(X)|\leq 64$ (or $b_2(X)=104$ and $sign(X)=-80$), then $X$ has non-vanishing 
refined invariants.
\end{theorem}

The special cases of vanishing first Betti numbers and $sign(X)=-16$, $sign(X)=-32$ and
$sign(X)=-48$ were obtained in \cite{MorganSzabo}, \cite{FKM} and \cite{FKMM}, respectively.

\subsection{Symplectic structures with $c_1=0$}

Let $X$ be a $K$-oriented  $4$-manifold, which is both symplectic and spin.  
This means the canonical $\sc$-structure coming with the symplectic structure
has vanishing integral first Chern class. For such a manifold one can combine the considerations
above with Taubes' result (\ref{TaubesWitten}). 

According to the Kodaira classification of complex surfaces \cite{BPV} there are only 
three families of complex surfaces with $c_1=0$. The obvious families are the
simply connected $K3$-surfaces and the tori. Furthermore, there are primary
Kodaira surfaces with first Betti number $3$. The first two families are K\"ahler,
hence symplectic by default. Other symplectic, non-K\"ahler and even non-complex 
manifolds with $c_1=0$ and Betti numbers
$2$ and $3$, were constructed \cite{Thurston}, \cite{FGG}, \cite{Geiges}.

\begin{theorem}{\rm (}\cite{MorganSzabo}{\rm)}
Let $X$ be a symplectic $4$-manifold with vanishing first Betti number and
with trivial canonical line bundle. Then $sign(X)=-16$.
\end{theorem}

{\bf Proof.} The vanishing of $c_1$ forces $X$ to be spin with $sign(X)=-16n$ and
$b^+=4n-1$. The monopole map is $Pin(2)$-equivariant, and after moding out the $\T$-action
as in \ref{proj}, we obtain a stable $\Z/2$-equivariant map in $\{P(\C^{2n}), S(V^{4n-1})\}^{\Z/2}$.
The action is free on both spaces and the spaces are of the same dimension. 
This allows to apply
a $\Z/2$-equivariant version of the Hopf theorem \cite{tD}, p 126.
According to this theorem, the (nonequivariant) degree of such a map is determined modulo $2$.  
As we will see, this degree, which is the Seiberg-Witten invariant, can be an odd number 
only in the case $sign(X)=-16$. Taubes' theorem \ref{TaubesWitten} then completes the argument.
To show that the degree is even for $n>1$, it suffices because of Hopf's theorem to 
exhibit an element in  $\{P(\C^{2n}), S(V^{4n-1})\}^{\Z/2}$ which has even degree. Here it is:
The $n$-th power $\eta^n$ of the $Pin(2)$-equivariant Hopf map 
induces a $\Z/2$-equivariant map $P(\C^{2n})\to S(V^{3n})$. Composed with the
inclusion $S(V^{3n})\to S(V^{4n-1})$ we get a map of degree zero for $n>1$.\qed

Here is an immediate corollary:

\begin{corollary}
A symplectic $4$-manifold with finite fundamental group
and with trivial canonical line bundle  is homeomorphic
to a $K3$-surface.
\end{corollary}

It is well-known that there are infinitely many different smooth structures
on the topological $4$-manifold underlying a $K3$-surface. 
The infinitely many smooth structures which come from complex 
K\"ahler structures were classified in \cite{Bauer92}. 
There are infinitely many more smooth structures which come from 
symplectic ones, compare \cite{GompfStipsicz}, p. 396f. And there are again
infinitely many more smooth structures which don't allow for a symplectic
structure at all, compare \cite{FintushelStern}. Nevertheless, amongst all
these smooth structures only the $K3$-surface seems to be 
known to carry a symplectic structure with $c_1=0$. 
The analogy to the Kodaira-dimension zero case in the Kodaira-classification
therefore is tantalizing:

\begin{question} Are symplectic $4$-manifolds with $c_1=0$ necessarily
either parallelizable or $K3$-surfaces? 
\end{question}

\subsection{Group actions}\label{Markus}

The stable cohomotopy approach 
does not rely on transversality results and therefore
seems suitable for considering
group actions on $4$-manifolds. 
In the discussion below, which closely follows \cite{Szymik}, 
the first Betti number of the manifolds will always 
be zero.  A compact Lie group acting on a $4$-manifold $X$ is supposed to 
preserve its $K$-orientation.

\begin{theorem}{\rm (}\cite{Szymik}{\rm)}
Let $G$ act on the $K$-oriented $4$-manifold $X$.  
There is a central extension $\mathbb G$ of the group
$G$ by the torus $\T$, such that the monopole map $\mu:\CA\to\CC$ is 
$\mathbb G$-equivariant. The associated element
$[\mu]_G\in \pi^0_{{\mathbb G}, \CU}({\mathrm {ind}\,l})$
restricts to the $\T$-equivariant stable cohomotopy
invariant.
\end{theorem}

When considering free actions of finite groups, one gets into a situation which very
much resembles Galois theory. The quotients $X/H$
by the various subgroups $H<G$ are $4$-manifolds carrying a residual action
of the Weyl group $WH=N_GH/H$.

\begin{theorem}{\rm (}\cite{Szymik}{\rm)}
Let $X$ be a $K$-oriented $4$-manifold with a free action of a finite group $G$ and $H<G$ a 
subgroup.
The set $J(H,\fs_X)$ of $\sc$-structures on the quotient $X/H$ which pull back to the
given $\sc$-structure $\fs_X$ on $X$ can be canonically identified with the
set of subgroups of $\mathbb G$ which map isomorphically to $H$ under the projection
to $G$. For $j\in J(H, \fs_X)$ the invariant
$[\mu(X/H,\fs_j)]_{WH}$ can be identified with the restriction of
$[\mu(X,\fs)]_G$ to the fixed points of $H(j)<{\mathbb G}$.
\end{theorem}

In particular, stable cohomotopy invariants
of oriented $4$-manifolds with finite fundamental group are determined by 
equivariant stable cohomotopy invariants of simply connected $4$-manifolds.

One can combine all the 
restrictions to fixed points into a comparison map
$$
\pi^0_{\mathbb G, \CU}({\mathrm {ind}\,l}(X,\fs))\to \bigoplus_{(H)\leq G}\,\bigoplus_{J(H,\fs_X)} 
H^0(WH; \pi^0_{\T,\,\CU^{H(j)}}({\mathrm {ind}\,l}(X/H; \fs_j))).
$$

Under certain conditions 
a general splitting result in equivariant homotopy theory implies that this
comparison map is an isomorphism after localisation away from the order of the group.
This splitting theorem can be applied for example if both   
$b^+(X/H)>1$ holds for any subgroup of $G$ and the index of the Dirac operator can 
be represented by an actual representation. 
So in this case kernel and cokernel are torsion groups with nonzero $p$-primary parts only for
those primes which do divide the order of $G$.

Finally, let's restrict to the case of a group of prime order $p$. Again the 
case where the $K$-orientation on $X$ comes from an almost complex structure
is easy to handle.

\begin{theorem}{\rm (}\cite{Szymik}{\rm)}
If the group $G$ of prime order $p$ acts freely on the almost complex manifold $X$,
then the invariant $[\mu(X,\fs)]_G$ is completely determined by the non-equivariant
invariants for $X$ and for $X/G$. Among the latter, the relation
$$[\mu(X,\fs)]\equiv \sum_{J(G,\fs)}[\mu(X/G, \fs_j)]\, mod\,p$$
is satisfied. 
\end{theorem}

The comparison map, however, 
is not injective in general. This is proved in \cite{Szymik}
using Adams spectral sequence calculations.
To find geometrical applications for these homotopy theoretical
computations looks like a challenging problem.


\section{Final remarks}

There is no chance to determine stable cohomotopy invariants by direct computation.
This seems obvious. So the only way to get further information out of the
monopole map is through a better conceptual understanding. 

Any improvement in our knowledge about the 
groups which arise as equivariant stable cohomotopy groups in this field
could help as a guideline to computing invariants  
as well as to constructing $4$-manifolds. 
We know disturbingly few examples of non-vanishing refined invariants. 
All examples known at the moment are powers of the Hopf map $\eta$.
Actually this reflects the fact that $\eta$
by Pontrjiagin-Thom describes the Lie group framing of the group $\T$ acting.

A hypothetical way to realizing other stable cohomotopy elements 
was pointed out in \ref{Birgit}: Construct a minimal counterexample to the
$\frac{11}8$-conjecture! Now we know, where to start the search (\ref{minimalspin}). 
It doesn't look like that hopeless an enterprise anymore.

It were symmetry considerations which lead to \ref{Birgit}. Indeed,
symmetry considerations may be a key to further progress. 
Let's dwell upon it a little more.
One can consider the monopole map as a map between infinite dimensional
bundles over some configuration space $Conf(X)$ consisting of all the choices made: metrics,
$\sc$-connections, harmonic $1$-forms. 
There is a symmetry group $\mathbb G$ acting: It is an extension
of the subgroup of the diffeomorphism group preserving the $K$-orientation
by some gauge group. Ideally, the monopole map can be understood as an Euler class of the
virtual index bundle
in a ``proper stable $\mathbb G$-equivariant cohomotopy group'' 
$$\pi^0_{\mathbb G}(Conf(X); {\mathrm{ind}\,l})$$
with twisting in
an element of ``proper $\mathbb G$-equivariant $KO$-theory''. 
The space $Conf(X)$ is the classifying space for proper $\mathbb G$-actions.
The obvious map $E{\mathbb G}\to Conf(X)$ from the classifying space of free actions
induces a ``Segal map''
$$\pi^0_{\mathbb G}(Conf(X); {\mathrm{ind}\,l})\to \pi^0_{\mathbb G}(E{\mathbb G}; {\mathrm{ind}\,l}).
$$
In analogy to the compact Lie group case one would expect
the latter group to be isomorphic (or at least related)
to non-equivariant stable cohomotopy
$\pi^0(B{\mathbb G};  {\mathrm{ind}\,l})$.
Now the classifying space $B{\mathbb G}$ of the group $\mathbb G$ 
indeed classifies parametrized
families of $K$-oriented $4$-manifolds. 
The image of the monopole class in this last group therefore is the universal
parametrized stable cohomotopy invariant. 
Of course, everything here is ill defined and probably cannot be made precise at all.
However, it can be made precise for compact approximations, i.e. for compact subgroups
of $\mathbb G$ or for finite equivariant subcomplexes of $Conf(X)$. This might lead 
to information on the diffeomorphism groups
of $4$-manifolds. Already the case of the four-dimensional sphere looks
interesting. 

Interesting first results in this direction, relating diffeomorphisms of $4$-manifolds to
parametrized stable cohomotopy invariants over the $1$-sphere,  can be found in a recent preprint \cite{MingXu}.

The space $Conf(X)$ might also be of interest for considering the behaviour of the stable cohomotopy
invariants at its ``boundary'', i.e. study the behaviour
of the maps under degeneration of the manifolds.

Another challenging direction of research is to find homotopy interpretations of Donaldson invariants 
and of Gromov-Witten invariants and to relate these concepts. At the moment these seem to be totally
out of reach.

Most urgently needed, however, are more general concepts  of gluing. 
Ideally, there should be relative invariants for manifolds with boundaries defining a 
``stable homotopy'' quantum field theory.  
Sadly enough, such concepts are still missing.

A number of speculative preprints on this topic have been circulating. Just at the time 
of writing at least one of them is being published.
The problems (cf. \ref{spectraI}, \ref{htop}) concerning well-definedness 
had been pointed out repeatedly to the authors as well as to the editors.

\bigskip{\it Acknowledgement:} I am grateful to M. Szymik for critical and helpful
comments.

\end{document}